%% file: Disc_payoff_21-07_1500.tex
\documentclass[review]{siamonline0516}

\input{title_page.tex}

\ifpdf
\hypersetup{
  pdftitle={\TheTitle},
  pdfauthor={\TheAuthors}
}
\fi

\newsiamremark{remark}{Remark} 
\crefname{remark}{Remark}{Remark} 

\newtheorem{prop}[theorem]{Proposition}
\usepackage{dsfont}

\begin{document}

\maketitle

\begin{abstract}
  We consider the   pricing problem related to payoffs that can have discontinuities of polynomial growth. The asset price dynamic is modeled  within the  Black and Scholes framework characterized by a stochastic volatility term driven by a fractional Ornstein-Uhlenbeck process.
  In order to solve the aforementioned problem, we consider three approaches.
 The first one consists in a suitable transformation of the initial value of the asset price, in order to eliminate possible discontinuities. Then we discretize both the Wiener process and the fractional Brownian motion and estimate the rate of convergence of the related discretized price to its
real value, the latter one being impossible to be evaluated analytically.
The second approach consists in considering   the conditional expectation with respect to the entire trajectory of the fractional Brownian motion (fBm). Then we derive a closed  formula which involves only  integral functional depending on the fBm trajectory, to evaluate the price; finally we   discretize the fBm and estimate the  rate of convergence of the associated numerical scheme to the option price.  In both cases the rate of convergence is the same and equals $n^{-rH}$, where $n$ is a number of the points of discretization, $H$ is the Hurst index of fBm,  and $r$ is the H\"{o}lder exponent of volatility function. The third method consists in calculating the density of the integral functional depending on the trajectory of the fBm via Malliavin calculus also providing the option price in terms of the associated probability density.

\end{abstract}



\newcommand{\s}{\mathbf{S}}
\newcommand{\Y}{\mathbf{Y}}
\newcommand{\X}{\mathbf{X}}
\newcommand{\W}{\mathbf{W}^1}
\newcommand{\WW}{\mathbf{W}^2}
\newcommand{\R}{\mathbb{R}}
\newcommand{\N}{\mathbb{N}}
\newcommand{\E}{\mathrm{E}}

\newcommand{\I}{\mathds{1}}

 \section{Introduction} Starting with the pioneering works by Hull and White \cite{HW87}   and     Heston \cite{heston}, models of financial market whose asset prices include stochastic volatility  have been the subject of an intensive research activity, which is still vibrant    from analytical, computational and statistical points of view.    Of course, option pricing is one of most relevant  problems.
In the latter context, stochastic volatility models are widely used because of their flexibility.
  Concerning the question  how to model stochastic volatility,   note that there are approaches in terms of   Gaussian (\cite{nikolato}, \cite{zhu}), non-Gaussian (\cite{barn1}, \cite{barn} ), jump-diffusion and L\`{e}vy processes (\cite{kypr}, \cite{tank}), as well as time series (\cite{carras}, \cite{palm}, \cite{shep}).  Our references are not in any way intended to be exhaustive or complete, we only illustrate the availability of different models. We would also like to mention the    books \cite{fouque},  \cite{kahel},  \cite{sat} and references therein, as well as  the paper \cite{AN15} which in some sense was a starting point for our considerations.  Furthermore, the models of financial market where the asset price includes stochastic volatility with long
memory in the volatility process is a subject of extensive research acrivity, see, e.g., \cite{boller}, where a wide  class of fractionally integrated GARCH and EGARCH models for characterizing financial market volatility was studied, \cite{comte} for affine fractional stochastic volatility models and \cite{viens}, where the Heston model with fractional Ornstein--Uhlenbeck  stochastic volatility was studied. As   was mentioned in \cite{comte}, long memory included into the volatility model  allows   to explain some option pricing puzzles such as steep volatility smiles in long term options and co-movements between implied and realized
volatility.  Note also that as a rule, the option pricing in stochastic volatility models needs some approximation procedures including Monte--Carlo methods.

 The present paper   contains a comprehensive and diverse approach to the exact and approximate option pricing of the asset price model that is described by the linear model with stochastic volatility, and volatility is driven by fractional Ornstein--Uhlenbeck process with Hurst index $H>\frac{1}{2}$. For technical simplicity we assume that the Wiener process driving the asset price and the fractional Brownian motion driving stochastic volatility are independent.  In these features, our model is similar to the model considered in \cite{viens}. However, the significant novelty of our approach  is that we consider three possible levels of  representation and approximation of option price, with the corresponding rate of convergence of discretized option price to the original one. Another novelty is that we can rigorously treat the class of discontinuous
  payoff functions of polynomial growth. As an example, our model allows to analyze linear combinations of digital and call    options. Moreover, we provide, for the first time in literature to the best of
  our knowledge, rigorous estimates for the rates of convergence of option prices for   polynomial discontinuous payoffs
   $f$  and H\"{o}lder volatility coefficients, a crucial feature considering settings for which exact pricing is
  not possible. \\

 The first level    corresponds to the case when    the price is presented as the functional of both driving stochastic processes, the Wiener process and the fBm,  and we dicretize and simulate both the trajectories of the Wiener process and of the fBm (double discretization) and estimate the rate of convergence for the discretized model. In these settings we apply the elements of the Malliavin calculus, following \cite{AN15}, to transform the option price to the  form that  does not contain discontinuous functions.     The second level corresponds to the case when we discretize and simulate only  the trajectories of the  fBm involved in  Ornstein--Uhlenbeck  stochastic volatility process (single discretization), basically  conditioning on the stochastic volatility process, then calculating the corresponding option price  as a functional of the trajectory of the fBm, and finally   estimating the rate of convergence of the discretized price.  This approach allows to  simulate only the trajectories of the fBm.   Corresponding simulations are presented and compared to those  obtained by the first   level. We conjecture  that the single discretization  gives better simulation results.  The third level potentially permits to avoid simulations, because it is possible to provide an analytical expression  for the option price, as an integral including  the density  of the functional which depends on stochastic volatility. Nevertheless, the density   whose existence we can prove in the framework of Malliavin calculus, is rather complicated from the computational point of view, therefore this   level  is more of theoretical nature.

  Taking into account previously mentioned approaches and techniques, our contribution is concerned with the treatment of a financial market, characterized by a finite maturity time $T$, and composed by a risk free bond, or bank account, $\beta=\{ \beta_t , t \in [0,T]\}$, whose dynamic reads as  $\beta_t=e^{\rho t}$,
where $\rho \in \mathbb{R}^+$ represents the  risk free interest rate, and a risky asset $S=\{ S_t, t \in [0,T]\}$ whose stochastic price  dynamic is defined, over the probability space  $\{\Omega, \mathcal{F},\mathbb{F}=\{\mathcal{ F}_t\}_{t\in [0,T]},\mathrm{P}\}$, by the following system of stochastic differential equations
 \begin{align}
  dS _t & = b S  _t dt + \sigma (Y _t) S_t dW  _t, \ \ \
  \ \ \ \ \ \ \ \ \  \label{chafe}
  \\
  dY _t & = - \alpha Y _t dt + d B ^H _t,\;   t \in [0,T].   \label{snide}
 \end{align}
Here $W=  \{ W_t, t \in [0,T]\}$ is a standard Wiener process, $b \in \R, \alpha \in \R^+$, are constants, while $Y=  \{ Y_t,t \in [0,T]\}$  characterizes  the stochastic volatility term of our model, being the argument of the function $\sigma$.
The process $Y$ is Ornstein-Uhlenbeck, driven by  a fractional Brownian motion
  $B^H= \{ B^H_t, t \in [0,T]\}$,
of Hurst parameter $H > \frac 12$, assumed to be independent of $W$.
Let us recall that    fBm  is a centered Gaussian process with covariance function $\E B ^H _tB ^H _s=\frac{1}{2}(s^{2H}+t^{2H}-|t-s|^{2H}).$ Moreover, due to the Kolmogorov theorem,  any  fBm has a modification such that its trajectories are almost surely H\"{o}lder function  up to order $H$. In what  follows we shall consider such modification. Moreover, it is well-known that for $H > \frac 12$, fBm has a long memory. This is   suitable for  stochastic volatility which represents the memory of the  model. We would also like to remind that a market model as   one described by the system of equations   \cref{chafe}, \cref{snide}, is incomplete   because of  two  sources of uncertainty, wether or not it is  arbitrage-free.

Therefore, in what follows we focus our attention on the   so called physical, or {\it real world}, measure, instead of using an equivalent martingale one.
Note, however, that in the case when the market is indeed arbitrage-free and there exists a minimal martingale measure,  the stock prices evaluated  w.r.t. the minimal martingale measure, resp. w.r.t. the objective measure,  differ only  to the non-random coefficient $e^{(b-\rho)t}$, as it happens in the standard Black--Scholes framework. For the discussion of conditions for the absence of arbitrage in the markets with stochastic volatility    see, e.g., \cite{kuch}.
As concerns the payoff function, we consider a measurable one defined by $f:\R^+\rightarrow \R^+$, and depending on   the value $S_T$ of the stock at maturity time $T$. It follows that our main goal is to calculate and approximate $\mathrm{E}f(S_T)$ with respect to the   different aforementioned   levels, also providing  rigorous estimates for the   corresponding rate of convergence for the first and second levels.

 The paper is organized as follows: in   \cref{results}
we give additional assumptions on the components of the model and formulate auxiliary results;   \cref{malliavin} contains the necessary elements of the Malliavin calculus that will be  used later;   \cref{1st level} contains the main results on the rate of convergence of the discretized option pricing approach, when we simulate  the trajectories of the Wiener process and  of the fBm;  \cref{2nd  level} contains the main results concerning the rate of convergence of the discretized option pricing problem when conditioning on  the trajectories of  the fBm, hence only simulating its trajectories;  \cref{3rd level} is devoted to the analytical derivation of the option price in terms of the density of the volatility functional, without trajectories simulations; the proofs    are collected in    \cref{proofs}; finally,   \cref{sims} provides the computer simulations associated to the approaches described in    \cref{1st level} and \cref{2nd  level}.

 \section{Model of asset price and payoff function: additional assumptions, auxiliary properties}\label{results}

Throughout the paper we assume  that payoff function $f:\R^+\rightarrow \R^+$ satisfies the following conditions:\\
 $(\mathbf{A})$ \begin{itemize}
\item[]$(i)$
 $f$ is a measurable
 function of polynomial growth,
 \[
 f(x) \leq C _f (1+ x^p), \ \ \ x \geq 0,
 \]
for some constants $C_f>0$ and $ p > 0$.
\item[]$(ii)$ Function $f$ is locally Riemann integrable, possibly, having discontinuities of the first kind.
\end{itemize}
Moreover we assume that the function   $\sigma:\R\rightarrow \R$  satisfies the following conditions:\\
$(\mathbf{B})$  there exists $C_ \sigma >0$ such that
\begin{itemize}
\item[]$(i)$ $\sigma$ is
  bounded away from $0$, $\sigma (x) \geq \sigma _{\min}>0$;
  \item[]$(ii)$ $\sigma$ has moderate polynomial growth, i.e., there exists
 $q   \in (0,1)$
     such that
 $$\sigma (x)
   \leq C_ \sigma (1 + |x|^{q}),\;x\in\R;$$
   \item[]$(iii)$   $\sigma$    is uniformly H\"older continuous, so that there exists $r\in (0,1]$
     such that
   $$|\sigma (x)- \sigma (y)| \leq C_\sigma |x-y|^{r},\;x,y\in\R;$$

\item[]$(iv)$   $\sigma\in C(\R)$   is differentiable a.e. w.r.t. the Lebesgue measure on $\R$, and its derivative is of polynomial growth: there exists
 $q'>0$   such that
 $$|\sigma' (x)|
   \leq C_ \sigma (1 + |x|^{q'}), $$ a.e. w.r.t. the Lebesgue measure on $\R$.
\end{itemize}
\begin{remark} 1) Concerning the relations between properties $(ii)$ and $(iii)$, note that we allow $r =1$ in $(iii)$ whereas
	$(ii)$ follows from $(iii)$ only in the case  $r<1$.

 2) Concerning the relations between properties $(iii)$ and $(iv)$, neither of these properties implies the other one unless  $r=1$. Indeed, on the one hand, a typical trajectory of a Wiener process is H\"{o}lder up to order $\frac{1}{2}$ but nowhere differentiable, on the other hand,  even continuous differentiability does not imply the uniform H\"{o}lder property.
\end{remark}

According to \cite{nVV99}, fBm admits a compact interval representation via some Wiener process $B$,  specifically,
  \begin{equation}\label{eq:hoe}
   B^H _t = \int\limits _0 ^t k(t,s)  d B_s,\;\;  k (t,s) =
   c _H
   s^{\frac 12 - H}
   \int\limits _s ^t u^{H - \frac 12} (u-s)^{H -\frac 32} du\I_{s<t},
  \end{equation}
 with
 $c _H = (H - \frac 12 )\left( \frac{2H \Gamma (\frac 32 - H)}
   {\Gamma (H+\frac 12) \Gamma (2 - 2H)} \right)^{1/2}$.
Obviously, the  processes $B$ and $W$ are independent. The next result is almost evident, however, we formulate it and even give a short proof for the reader's convenience.

\begin{lemma}\label{lem:lemma 2.1}
\begin{itemize}
\item[] $(i)$  Equation \cref{snide}   has a unique solution of the form
\begin{equation}\label{eq:solution1}  Y_t=Y_0e^{-\alpha t}+\int_0^te^{-\alpha(t-s)}dB_s^H.\end{equation}
Moreover, for any $\alpha>0$ and any $\beta<2$ \begin{equation}\label{eq:exp-mom-sup}\mathrm{E}\exp\{\alpha\sup_{t\in[0,T]}|Y_t|^\beta\}<\infty.\end{equation}
\item[] $(ii)$  Equation  \cref{chafe} has a unique solution of the form
\begin{equation}\label{eq:solution2} S_t=S_0\exp\left\{bt+\int_0^t \sigma(Y_s)dW_s-\frac{1}{2}\int_0^t \sigma^2(Y_s)ds\right\}.\end{equation}
Moreover, for any $m\in \mathbb{Z}$ we have  $\mathrm{E} (S_T)^m<\infty$, and  for any $m>0$ it holds   $\mathrm{E }(f(S_T))^m<\infty.$
\end{itemize}
\end{lemma}
\begin{remark}\label{polyn growth} We can generalize the last conclusion of   \cref{lem:lemma 2.1} to the following one: for any function $\psi=\psi(x):\R\rightarrow\R$ of polynomial growth $\sup_{t\in[0,T]}\mathrm{E }(|\psi(S_t))|<\infty.$
\end{remark}

 \section{Elements of Malliavin calculus and application to option pricing}\label{malliavin}

In what follows, we recall some basic definitions and results about Malliavin calculus, doing that we mainly refer to   \cite{Nua06}.
 Let   $W=\{W(t), t \in [0,T]\}$ be a Wiener process on the standard probability space  $\{\Omega, \mathcal{F},$  $\mathrm{F}=\{\mathcal{F}_t^{W}\}, t \in [0,T],
\mathrm{P}\},$ where  $\Omega=C([0,T], \mathbb{R}).$
Denote by $\widehat{C}^{\infty} \mathbb(R)$ the set of all infinitely differentiable functions with the derivatives of polynomial growth at infinity.

\begin{definition}
Random variables $\xi$ of the form $\xi=h(W(t_1), \ldots , W(t_n))$,
$$h=h(x^1,\ldots,x^n) \in \widehat{C}^{\infty}(\mathbb{R}^n),\;
 t_1, \ldots
t_n \in [0,T], \; n\geq 1 $$ are called smooth.
Denote by   $\mathcal{S}$ the class of smooth random variables.
\end{definition}

\begin{definition}
Let $\xi \in \mathcal{S}$. The stochastic derivative of
  $\xi$ at    $t$ is the random variable
\begin{equation*}
D_t\xi= \sum_{i=1}^n \frac{\partial h}{\partial x^i}(W(t_1), \ldots , W(t_n)) \I_{t\in [0,t_i]}, \quad t\in[0,T].
\end{equation*}
\end{definition}

Considered as an operator from $L^2(\Omega)$
to $L^2(\Omega ; L^2[0,T])$,
$D$ is a closable operator.
We use the same notation $D$ for its closure.
$D$ is known as the Malliavin derivative,
or the stochastic derivative.
The domain of the operator of the stochastic derivative
 is a Hilbert space $\mathrm{D}^{1,2}$
of random variables,   on which the inner product
(which coincides with the operator norm)
is given by
$$
\langle \xi, \eta\rangle _{1,2}=\mathrm{E}(\xi\eta)+\mathrm{E}(\langle D\xi, D\eta \rangle _H), \quad H=L^2([0,T],\R).
$$

Thus, the operator of stochastic derivative  $D$ is closed,
unbounded and defined on
a dense subset of the space  $L^2(\Omega)$ (see \cite{Nua06}).
The following statement is known as the chain rule.

\begin{prop} (\cite[Proposition 1.2.3]{Nua06}).
Let $\varphi: \R ^m \to \R$
be a continuously differentiable function with
bounded partial derivatives. Suppose that
$\xi = (\xi_1,...,\xi_m)$ is a random vector whose components
belong to $\mathrm{D}^{1,2}$. Then
$\varphi (\xi) \in \mathrm{D}^{1,2}$ and
\[
 D\varphi (\xi)  =
 \sum\limits _{i= 1} ^m \partial _i \varphi (\xi) D \xi_i.
\]
\end{prop}

Denote by $\delta$ the operator adjoint to  $D$ and by  $\text{Dom}\; \delta $ its domain. The operator $\delta$ is unbounded in $H$ with   values in  $L^2(\Omega)$ and such that
\begin{itemize}
\item[]$(i)$ $\text{Dom}\; \delta $ consists of square-integrable random processes   $u\in H,$ satisfying
$$
\left|\mathrm{E}\left(\left\langle D\xi, u\right\rangle_H \right)\right|\le C(\mathrm{E}( \xi^2))^{1/2},
$$
for any $\xi\in\mathrm{D}^{1,2}$, where $C$ is a constant depending on  $u$;
\item[]$(ii)$ If $u$ belongs to  $\text{Dom}\; \delta $, then   $\delta(u)$ is an element of  $L^2(\Omega)$ and
$$
\mathrm{E}\left(\xi\delta(u)\right)=\mathrm{E}\left( \left\langle D\xi, u\right\rangle_H\right)
$$
for any $\xi\in\mathrm{D}^{1,2}$.
\end{itemize}
 The operator $\delta$ is closed.  Consider the space   $L^{1,2}=L^2([0,T],\mathrm{D}^{1,2})$  with the norm $||\cdot||_{L^{1,2}},$ where
\begin{equation*}
||u||^2_{L^{1,2}}=\mathrm{E}\left(\int_0^T u^2_tdt+\int_0^T\int_0^T \left(D_s u_t \right)^2 dtds\right).
\end{equation*}
If $u \in L^{1,2},$ then the integral $\delta(u)$ is correctly defined and  $$\mathrm{E}\left(\int_0^T u_t dW_t\right)^2\leq ||u||^2_{L^{1,2}}$$ (see
\cite{Nua06}). In this case operator $\delta(u)$ is called the Skorokhod integral of the process $u$ and is denoted by
$$\delta(u)=\int\limits_0^T u_t dW_t.$$
To apply Malliavin calculus to the asset price $S$, note that we have a two-dimensional case with two independent Wiener processes $(W, B)$. With evident modifications, denote by $(D^W, D^B)$ the stochastic derivative
  with respect to the two-dimensional Wiener process  $(W, B)$.
Denote also $$X (t)  = \log S (t) =
\log S _0 + bt -
\frac 12 \int\limits _0 ^t \sigma ^2 (Y _s)ds
 + \int\limits _0 ^t \sigma  (Y _s)dW_s.$$
 \begin{lemma}\label{lem:scrounge}
\begin{itemize} \item[]$(i)$ The stochastic derivatives of the fBm $B^H$ equal to
 $$D^W_u B^H _t  = 0, \;\; D^B_u B^H _t  = k(t,u).$$
\item[]$(ii)$ The stochastic derivatives of $Y$ equal to
 \begin{equation} \label{eq:vertigo}
\begin{split}
D^W_u Y _t  = 0, \quad
 D^B_u Y _t
 =
  c _H e^{-\alpha t} u ^{1/2 - H}
  \int\limits _u ^ t e^{\alpha s} s ^{H - 1/2} (s-u) ^{H - 3/2}ds  \I_{u <t}.
  \end{split}
\end{equation}
\item[]$(iii)$ The stochastic derivatives of $X$ equal to
\begin{align}\label{libelous}
     D _u^W X _t   = \sigma (Y _u)\I_{u <t}, \;   D _u^B X _t =   \left(
   -  \int\limits _0 ^t  \sigma (Y _s) \sigma ' (Y _s)
   D _u ^B Y _s ds +
   \int\limits _0 ^t   \sigma ' (Y _s)
   D _u ^B Y _s dW_s
    \right)  \I_{u <t}.
 \end{align}
\end{itemize}
\end{lemma}
 \begin{lemma}\label{lem:slander}
   The laws of $S_T$ and $X_T$
   are absolutely continuous with
   respect to the Lebesgue measure.
  \end{lemma}
  From now on, we denote $C$ any constant whose value is not important and can change from line to line and even inside the same line. Throughout the paper,   $C$ cannot depend on $n, t, s$, but can depend on $\sigma, H, T, Y_0, S_0, \alpha, b, p, r, q, q',f$ and other parameters specified in the problem.
 In what follows we need the statement contained in the next remark.
 \begin{remark} \label{important rem} The chain rule of stochastic differentiation can be extended to the wider class of functions in the following way. Applying Proposition 1.2.4 from \cite{Nua06} and the related remark, we get  that in the case when the function $\varphi$ is Lipschitz  and has a derivative a.e. w.r.t. the Lebesgue measure on $\R$, and the law of r.v. $\xi$ is absolutely continuous with respect
to the Lebesgue measure on $\R$, then $\varphi(\xi)$ has a stochastic derivative and $D\varphi(\xi)=\varphi'(\xi)D\xi$ a.e. w.r.t. the Lebesgue measure on $\R$.
 Now, consider the   stochastic differentiation  of the functions of Ornstein--Uhlenbeck process $Y$. Let $\varphi$ be locally Lipschitz, with
 both  and $\varphi$ and $\varphi'$ being of polynomial growth.
 Then
 $\varphi_n(x)=\varphi(x)\I_{ |x|\leq n }+
 (\varphi(-n))\I_{ x<- n}+\varphi(n)\I_{  x>n }$
 is Lipschitz, has a derivative $\varphi _n '$ a.e. w.r.t. the Lebesgue measure, and moreover,
  there   exists a polynomial  $\bar \varphi$
  with non-negative coefficients such that
  $$|\varphi_n(x)|+|\varphi_n '(x)|\leq \bar \varphi(|x|), \;x\in\R,$$
  which implies that
\[
\mathrm{ E}|\varphi (Y_s) - \varphi _n (Y_s)|^2
 \leq 4 \mathrm{E }\left(\bar \varphi ^2 (|Y_s|) \I_{Y_s \notin [-n,n]}\right)
  \to 0.
\]
Furthermore, it easily follows from  \cref{eq:vertigo} that in fact $D^B Y_s$ is in $L^2([0,T])$.  Indeed,
\begin{equation}\label{eq:mime}
    0\leq D^B_u Y_s \leq C u^{1/2 - H } (s-u)^{H -1/2}\I_{u<s}.
  \end{equation}   Furthermore, $$\mathrm{E} \left(\max\limits _{ s \in [0,T]}\bar \varphi^2(|Y_s|)\right)<\infty,$$ due to the fact that  $\max\limits _{ s \in [0,T]}|Y_s|$ has exponential moments. Therefore,
\begin{equation*}\begin{gathered}
 \mathrm{E} \left(\int\limits _0 ^T
  ( \varphi ' (Y_s) D^B_u Y_s - D^B _u \varphi _n (Y_s))^2 du\right)\\
  \leq 4 \mathrm{E} \left(\max\limits _{ s \in [0,T]}\bar \varphi^2(|Y_s|)
  \I_{ \max\limits _{ s \in [0,T]} Y_s \notin [-n,n]}
  \int\limits _0 ^s  (D^B_u Y_s)^2 du
   \right) \to 0.
\end{gathered}\end{equation*}
 Previous results, together with the fact that $D$ is closed, imply that
 $D^B _u\varphi (Y_s) = \varphi ' (Y_s) D^B_u Y_s.$
 \end{remark}

Let us introduce the following notations:
 $g(y) = f(e^y)$,
 $F(x) = \int\limits _{0} ^x f(z)dz$
 and let
  $G(y) = \int\limits _{0} ^y g(z)dz$,
  $x \geq 0$, $ y \in \R$. Also, let
 \begin{equation}\label{eq:randvarZ}
 Z _T = \int\limits _0 ^T
   \sigma^{-1}(Y_u) dW_u.
  \end{equation}
Note that $Z_T$ is well defined because of condition  $(\mathbf{B})$, $(i)$. Now, analogously to \cite{AN15}, we   are in position to transform  the option price in such a way that it does not contain discontinuous functions.
    \begin{lemma}\label{lem:decry}
Under conditions $(\mathbf{A})$ and $(\mathbf{B})$ the
option price  $\mathrm{E} f(S_T)=\mathrm{E} g(X_T)$ can be represented as
\begin{equation}\label{eq:for f}
  \mathrm{E }f(S_T) = \mathrm{E }\left(\frac{F(S_T)}{S_T}\left( 1 + \frac{Z_T}{T}\right)\right).
 \end{equation}
 Alternatively,
 \begin{equation}\label{eq:for g}\mathrm{E} g(X_T) =\frac 1T \mathrm{E} \left(G(X_T)Z_T\right).
\end{equation}
\end{lemma}

\section{The rate of convergence of approximate  option prices in the case  when both Wiener process and fractional Brownian motions are discretized}\label{1st level}
In the present section we provide our first approach (first level) to the numerical approximation of the solution for the option pricing problem.
In particular, we are going to provide a double discretization procedure, with related simulations, with respect to both the Wiener process and the fBm, also estimating the rate of convergence for the corresponding approximated option prices to the real value given by    $\mathrm{E}f(S_T)$.

To pursue latter aim, let us introduce the following notation.
 For any  $n \in \N$ consider equidistant partition of the interval $[0,T]$: $t_i = t_i (n) = \frac{iT}{n}$,
 $i = 0,1,2,...,n$. Then we define the  discretizations of Wiener process $W$ and fractional Brownian motion $B^H$:

 \[
 \Delta W_i= W \left(t_{i+1}\right)-
  W \left(t_i\right),
  \]
  \[
  \Delta B^{H}_i= B^{H} \left(t_{i+1}\right)-
  B^{H} \left(t_i\right), i = 0,1,2,...,n.
 \]

Discretized  processes $Y$ and $X$,  corresponding to a given partition have the form
\[
  Y^n_{t_j}=
 Y_0 e^{-\alpha t_j} + e^{-\alpha t_{j-1}}
  \sum\limits _{i =0} ^{j-1}
  e^{\alpha t_i}
   \Delta B^H_i,
 \]
\[
 X^n_{ t_j}
 =
 X_0 + b t_j
 -\frac {1}{2n}
 \sum\limits _{i = 0} ^{j-1}
  \sigma ^2 \left(Y^n_{t_i}\right)
 +
 \sum\limits _{i = 0} ^{j-1}
  \sigma  \left(Y^n_{t_i}\right)
   \Delta W_{i}
\]
\[
 =  X _0 + b t_j  -
 \frac 12 \int\limits _0 ^{t_j} \sigma ^2 (Y^n _s)   ds
 + \int\limits _0 ^{t_j} \sigma (Y^n _s)    dW_s,
  \ \ \ j = 0,...,n,
\]
where  we put $ Y ^n _s =
Y^n_{t_i}$
for $s \in \big[ t_i, t_{i+1} \big)$.
Concerning the discretizaion of the term  $
 Z _T$  from \cref{eq:randvarZ}, it has a form
$Z ^n _T = \int\limits _0 ^T
   \frac{1}{\sigma (Y^n _s)   } dW_s.$
Eventually we define $S ^n _{t_j}=
  \exp \big\{ X ^n_{t_j}\big\}.$
 Three   lemmas below contain  all auxiliary bounds that are necessary in order to establish the main result.
\begin{lemma}\label{lem:disgruntle}\begin{itemize}
 \item[]$(i)$ For any $\theta>0$ there exists a  constant $C$ depending on $\theta$ such that for any  $s,t\in[0,T]  $

  \begin{equation*}
   \mathrm{ E} |Y_t - Y _s|^{\theta} \leq C \left| t-s \right|
    ^{\theta H}.
  \end{equation*}

 \item[]$(ii)$ For any $ \theta >0$ there exists a  constant $C$ depending on $\theta$ such that for any $0\leq j\leq n$
   \begin{equation*}
    \mathrm{E}\left|Y_{t_j} -
    Y^n_{t_j}\right|^{\theta} \leq
   Cn^{-\theta}.
 \end{equation*}
\item[]$(iii)$  For any $ \theta>0$ there exists a constant $C$ depending on $\theta$ such that for any $s\in[0,T]$
\begin{equation*}\begin{gathered}
\mathrm{E}|Y _s - Y ^n_{s} |^{\theta} =
\mathrm{E}|Y _s - Y ^n_{t_i} |^{\theta}
\leq C n^{-\theta H}.
\end{gathered}\end{equation*}

 $(iv)$ Approximating  process has uniformly bounded  moments:  for any $ \theta>0$
 \begin{equation}\label{eq:prop approxi}
    \sup_{s\in[0,T]}
    \mathrm{E} |Y^n _s|^\theta < \infty.
\end{equation}
\end{itemize}
\end{lemma}
\begin{remark}\label{boundinn} Using \cref{eq:prop approxi} and the fact that the approximating process $Y^n$ is Gaussian, we can prove similarly to   \cref{lem:lemma 2.1} and   \cref{polyn growth} that for any $m\in \mathbb{Z}$ $$\sup_{n\geq 1}\sup_{0\leq j\leq n}\mathrm{E}\left(S^n_{t_j}\right)^m<\infty.$$
\end{remark}
\begin{lemma}\label{lem:tattle}
 There exists a constant $C >0$  such that for any $n\geq 1$
 \begin{equation}\label{eq:haggle}
    \mathrm{E}(X_T - X ^n _T) ^2 \leq
   C n^{-2rH},
 \end{equation}
 and
 \begin{equation}\label{eq:bound for Z}
   \mathrm{ E}(Z_T - Z ^n _T) ^2 \leq
   C n^{-2rH}.
 \end{equation}

\end{lemma}
\begin{lemma} \label{lem:pull in horns}
 Under conditions $(\mathbf{A})$ and $(\mathbf{B})$
      we have the following upper bound: there exists a constant $C_F$ such that
      \begin{equation*}
       \mathrm{E} \left|\frac{F(S_T)}{S_T} -
       \frac{F(S _T ^n)}{S _T ^n} \right| ^2 \leq
       C_F\cdot n^{-2rH}.
      \end{equation*}

   \end{lemma}
Using previous lemmas, we are now in position to state the main result of this section, namely to provide the rate of convergence of discretized option prices to the exact one represented by $\mathrm{E }f(S_T)$, under double discretization.
 \begin{theorem}\label{thm:convergence thm} Let conditions $(\mathbf{A})$ and $(\mathbf{B})$ hold.
  There exists a constant $C$ not depending on $n$ such   that

       \begin{equation*}
  \left|\mathrm{E }f(S_T) - \mathrm{E }\left(
  \frac{F(S ^n _T)}{S ^n _T}
  \left(
  1 + \frac {Z ^n _T}{T}
  \right)
  \right) \right| \leq C n^{-rH}.
 \end{equation*}

   \end{theorem}
   \section{The rate of convergence of approximate  option prices in the case  when only  fractional Brownian motion is  discretized}\label{2nd level}
The present section is devoted to the implementation of the second approach (second level) to approximate the option price. It is based on the fact that in the case when $W$ and $B$ are independent, logarithm of asset price is conditionally Gaussian under the fixed trajectory of fractional Brownian motion. It allows to exclude Wiener process $W$ from the consideration and to calculate the option price explicitly in terms of the trajectory of fBm $B^H$. Respectively, we can discretize  and simulate only the trajectories of $B^H$ (single discretization).   \cref{thm:analytic expression thm} gives the explicit option pricing formula as the functional of the trajectory of fBm $B^H$, and   \cref{thm:second conv thm} gives the rate of convergence. Comparing to   \cref{thm:convergence thm},
we see that the rate of convergence admits the same bound, influenced by the behavior of volatility.

   Let us introduce the following notations:  let the covariance   matrix reads as follows
   \begin{equation*}
    C_{X,Z} = \begin{pmatrix}
    \sigma_{_Y} ^2 & T\\
    T & \sigma_{_Z} ^2
  \end{pmatrix},
   \end{equation*}
and let
   \begin{equation*}
     \sigma_{_Y} ^2 = \int\limits _0 ^T
    \sigma^2 (Y_s)ds, m_{_Y}=  X_0+bT- \frac 12 \sigma_{_Y} ^2, \sigma_{_Z} ^2 = \int\limits _0 ^T
    \sigma^{-2 }(Y_s)ds,  \Delta=|C_{X,Z}|=\sigma_{_Y} ^2\sigma_{_Z} ^2-T^2.
   \end{equation*}

Evidently, $\Delta\geq 0. $ We assume additionally that the following assumption is fulfilled.\\

$(\mathbf{C})$ $\Delta=\sigma_{_Y} ^2\sigma_{_Z} ^2-T^2>0$ with probability 1.\\

 Note that the random vector $$(X_T, Z_T)  =\left(X_0 + bT -
\frac 12 \int\limits _0 ^T \sigma ^2 (Y _s)ds + \int\limits _0 ^t \sigma  (Y _s)dW_s,\;\int\limits _0 ^T\sigma^{-1}(Y_s)  dW_s\right)$$
is  Gaussian conditionally to the given  trajectory $\{Y_t,t\in[0,T]\}$.
The conditional covariance matrix is $C_{X,Z}$.
Next lemma presents common   conditional  density of $(X_T, Z_T)$. Note that under assumption $(\mathbf{C})$ the distribution of $(X_T, Z_T)$ is non-degenerate in $\R^2$.
\begin{lemma}\label{lem:lem 5.1} Let assumption $(\mathbf{C})$ hold. Then the common conditional  density  $p_{X,Z}(x,z)$ of $(X_T, Z_T)$, conditionally to the  given trajectory $\{Y_t,t\in[0,T]\}$, equals
\begin{equation}\label{eq:boink}
    p_{X,Z}(x,z) =
    \frac{1}{2 \pi \Delta^{\frac{1}{2}}}
    \exp
    \Bigg\{
   - \frac{1 }{2
   \Delta
   }
   \left(
      \sigma _{_Z} ^2{(x-m_{_Y})^2}
    + \sigma_{_{Y}} ^2  {z^2}
    -
     {2T (x-m_{_Y})z }
   \right)
    \Bigg\}.
   \end{equation}
\end{lemma}
The next result states that option price can be presented as the functional  of $\sigma_{Y} ^2$  only.

   \begin{theorem}\label{thm:analytic expression thm}
 Under conditions $(\mathbf{A})$--$(\mathbf{C})$ the following equality holds:

\begin{equation}\begin{gathered}\label{eq:assuage}
 \E g(X _T)
 =
   (2\pi)^{-\frac{1}{2}}\int\limits _{\R } G(x)\mathrm{E}\left(
 \frac { (x - m _{_Y}) }{\sigma_{_Y} ^3}
        \exp  \Bigg\{
    -\frac{ (x - m _{_Y}) ^2}{2\sigma_{_{Y}} ^2}
       \Bigg\} \right)dx\\ =(2\pi)^{-\frac{1}{2}}\mathrm{E}\left((\sigma _{_Y})^{-1}\int\limits _{\R } G((x+m_Y)\sigma_{Y})xe^{-\frac{x^2}{2}}dx\right).
\end{gathered}\end{equation}

\end{theorem}
In order to state the main result of the present section, et us define the following quantities
\[
 \sigma_{_{Y,n}} = \int\limits _0 ^T
    \sigma^2 (Y ^n _s)ds,  \ \ \
    m_{_{Y,n}}=  X_0+bT- \frac 12 \sigma_{_{Y,n}} ^2.
\]

\begin{theorem}\label{thm:second conv thm}

Under conditions $(\mathbf{A})$, $(\mathbf{B})$,
and $(\mathbf{C})$  we have

\begin{equation*}
 \left| \E g(X _T) - (2\pi)^{-\frac{1}{2}}\int\limits _{\R } G(x)\mathrm{E}\left(
 \frac { (x - m _{_{Y,n}}) }{\sigma_{_{Y,n}} ^3}
        \exp  \Bigg\{
    -\frac{ (x - m _{_{Y,n}}) ^2}{2\sigma _{_{Y,n}} ^2}
       \Bigg\} \right)dx \right|
       \leq C  n  ^{-rH}.
\end{equation*}

\end{theorem}

\section{Option price in terms of density of the integrated stochastic volatility}\label{3rd level}
Applying    \cref{thm:analytic expression thm} and  equality \cref{eq:assuage}, we clearly see that the option price depends on the random variable  $\sigma_Y^2=\int\limits _0 ^T\sigma^2 (Y_s)ds$. Therefore it is natural to derive the  density  of this random variable. Since $\sigma_Y^2$  depends on the whole trajectory of the fBm $B^H$ on $[0,T]$, we apply Malliavin calculus in an attempt to find the density. First, establish some auxiliary results. For any $\varepsilon>0$ and $\delta>0$ introduce the stopping times $\tau_ \varepsilon =
 \inf\{t>0: |B^H_t| \geq \varepsilon \}$
 and $\nu_ \delta =
 \inf\{t>0: |{Y_t} - {Y_0}| \geq \delta \}$.
\begin{lemma}\label{lem:negative moments}
For any $l>0$   negative moments are well defined:    $\E(\nu_ \delta)^{-l}<\infty$.
 \end{lemma}

Now we introduce additional assumptions on the function $\sigma$.

$(\mathbf{D})$ The function $\sigma \in C^{(2)}(\R)$, its  derivative $\sigma'  $ is strictly nonnegative, $\sigma'(x)>0, \;x\in\R,$ and $\sigma'$, $\sigma''$ are of polynomial growth.
\begin{lemma} \label{lem:lurid} Under assumptions $(\mathbf{B})$ and $(\mathbf{D})$
 the stochastic  process
  \[
   \frac{D^B  \sigma ^2 _{_Y}}{||D^B  \sigma ^2 _{_Y}||^2 _{H}}=\left\{\frac{D^B_t  \sigma ^2 _{_Y}}{||D^B  \sigma ^2 _{_Y}||^2 _{H}},\;t\in[0,T]\right\}
  \]
 belongs to the domain $\text{Dom}\; \delta$ of the Skorokhod integral $\delta$.

 \end{lemma}
 Denote $\eta=(||D^B  \sigma ^2 _{_Y}||^2 _{H})^{-1},\; l(u,s)= c _He^{-\alpha s}
  \int\limits _u ^ s e^{\alpha v} v ^{H - 1/2} (v-u) ^{H - 3/2}dv $,
  $\kappa(y) = \sigma (y)\sigma'(y)$.
\begin{theorem}\label{thm:density fBm}\begin{itemize}
\item[]$(i)$ The density $p_{\sigma ^2 _{_Y}}$ of the random variable
$\sigma ^2 _{_Y}$ is bounded, continuous and given  by the following formulas

\begin{equation}\begin{gathered}\label{eq:loin-1}
  p_{\sigma ^2 _{_Y}} (u) = \E \left[ \I_{\sigma ^2 _{_Y} >u  }
  \delta \left(
  \frac{D^B \sigma  ^2 _{_Y}}{||D^B  \sigma ^2 _{_Y}||^2 _{H}}\right)
  \right],
 \end{gathered}\end{equation}
 where the Skorokhod integral is in fact reduced to a Wiener integral,  $$\delta \left(
  \frac{D^B \sigma  ^2 _{_Y}}{||D^B  \sigma ^2 _{_Y}||^2 _{H}}\right)=2\eta\int_0^T\kappa(Y _s)\left(\int_0^su ^{1/2 - H}l(u,s)dB_u\right)ds-\int_0^TD^B_u \eta  D^B_u (\sigma  ^2 _{_Y})du.$$

 \item[] $(ii)$
  The option price $\E g(X_T)$  can be represented as the integral with respect to the density $p_{\sigma ^2 _{_Y}} (u)$
  defined by \cref{eq:loin-1} as follows:
   \[
   \E g(X_T)
   \]
   \[=
    (2\pi)^{-\frac12}T \int\limits _{\R}  G(x)  \int_ {\R}
       \frac { (x + u/2 - X_0 - bT ) }{u ^3  }
        \exp  \Bigg\{
    -\frac{ (x + u/2 - X_0 - bT) ^2}{2u ^2}
       \Bigg\} p_{\sigma ^2 _{_Y}} (u) du.
   \]
 \end{itemize}
 \end{theorem}

\section{Proofs} \label{proofs}

\textbf{Proof of   \cref{lem:lemma 2.1}}.
$(i)$ The representation \cref{eq:solution1} for the fractional Ornstein-Uhlenbeck process $Y$ is well known, see, e.g., \cite{CKM03}. It is a continuous Gaussian process with $\sup_{t\in[0,T]}\mathrm{E}(Y_t)^2<\infty.$ Then the  boundedness of any exponential moments of the form \cref{eq:exp-mom-sup} follows from  \cite{fern} and \cite{Tal94}.\\
$(ii)$ To establish the representation \cref{eq:solution2} for
 $S$, we need only to prove that the  integrals $\int_0^t \sigma(Y_s)dW_s$ and $\int_0^t \sigma(Y_s)S_sdW_s$ are well defined, while the form of the representation is obvious. To what concerns $\int_0^t \sigma(Y_s)dW_s$, it follows from   \cref{eq:exp-mom-sup} and condition $(\mathbf{B})$, $(i)$ that  $\int_0^t  \mathrm{E}\sigma^2(Y_s)ds\leq  C\int_0^t  \mathrm{E}(1+|Y_s|^{2q})ds< \infty.$ Moreover, obviously, the moments of any order are finite: $\sup_{t\in[0,T]}\mathrm{E}\sigma^{2n}(Y_t)<\infty$.
 Furthermore, the conditional distribution of the integral $\int_0^t \sigma(Y_s)dW_s$ given the whole trajectory $Y=\{Y_t, t\in[0,T]\}$
 is Gaussian.
 Since we have $q<1$  together with exponential inequality \cref{eq:exp-mom-sup}, then  for any $n \in \mathbb{Z}$, the following holds
\begin{equation}\begin{gathered}\label{eq:exp-finity}\sup_{t\in[0,T]}\mathrm{E}S_t^n\leq C\sup_{t\in[0,T]}\mathrm{E}\exp\left\{n\int_0^t \sigma(Y_s)dW_s\right\}=C\sup_{t\in[0,T]}\mathrm{E}\exp\left\{\frac{n^2}{2}\int_0^t \sigma^2(Y_s)ds\right\}\\ \leq C\sup_{t\in[0,T]}\mathrm{E}\exp\left\{C_\sigma\frac{n^2}{2}\int_0^t (1+|Y_s|^{2q})ds\right\}\leq C\mathrm{E}\exp\left\{\frac{C_\sigma T n^2}{2}\sup_{s\in[0,T]}|Y_s|^{2q}\right\} <\infty.\end{gathered}
 \end{equation}
 In particular, $$\int_0^T  \mathrm{E}(\sigma^2(Y_s)S^2_s)ds\leq T \left(\sup_{t\in[0,T]}\mathrm{E}S_t^4\sup_{t\in[0,T]}\mathrm{E}\sigma^{4}(Y_t)\right)^{\frac{1}{2}}<\infty$$
 whence the proof follows.

\qed

\textbf{Proof of  \cref{lem:scrounge}}. Statement $(i)$ follows directly from the definition of stochastic derivative and from the fact that $B$ and $W$ are independent. Similarly, the first equality in \cref{eq:vertigo}
is obvious since $Y$ is independent of $W$. Furthermore, integrating by parts   \cref{eq:solution1} and taking into account representation \cref{eq:hoe}, we get following equalities
 \begin{equation}\label{eq:OU direct} Y_t=Y_0e^{-\alpha t}+B_t^H-\alpha e^{-\alpha t}\int_0^te^{\alpha s}B_s^Hds,\end{equation}
 whence $$D_u^BY_t=k(t,u)-\alpha e^{-\alpha t}\int_u^te^{\alpha s}k(s,u)ds,$$
 where the kernel $k$ was introduced in \cref{eq:hoe}.  Note that the derivative $k_s'$ of the kernel   $k$ equals $$k_s'(s,u)=c_Hu^{\frac{1}{2}-H}s^{H-\frac{1}{2}}(s-u)^{H-\frac{3}{2}}\I_{u <s}.$$ It is an integrable function, therefore we can integrate by parts once again and get that $$D_u^BY_t=  e^{-\alpha t}\int_u^te^{\alpha s}k_s'(s,u)ds=c_H e^{-\alpha t} u ^{\frac{1}{2} - H}
  \int\limits _u ^ t e^{\alpha s} s ^{H - \frac{1}{2}} (s-u) ^{H - \frac{3}{2}}ds  \I_{u <t}.$$
The first equation from \cref{libelous} follows from the definition of stochastic derivative. Further,  $$D _u^B X _t =\left(
   - \frac{1}{2} \int\limits _0 ^t  D_u^B(\sigma^2 (Y _s))ds +
   \int\limits _0 ^t   D_u^B(\sigma  (Y _s))dW_s
    \right)  \I_{u <t}.$$
     Note that \cref{important rem}, together with Gaussian distribution of $Y_s$, allows to apply the chain rule to the continuous and a.e. differentiable function $\sigma$ taken at point  $Y_s$. Moreover, the result  can be written in the standard form, so that   $D_u^B(\sigma  (Y _s))=\sigma'(Y _s)D_u^B(Y_s)$, and
   \begin{equation}\label{eq:derivative}D_u^B(\sigma^2(Y _s))=2\sigma(Y _s)\sigma'(Y _s)D_u^B(Y_s)\end{equation}
    a.e. w.r.t. the Lebesgue measure on $\R$. Besides, similarly to proof of   \cref{lem:lemma 2.1}, we can apply properties $(\mathbf{B})$, $(ii)$ and $(iv)$, which, together with the  upper bound \cref{eq:mime}
   ensure that the integrals in \cref{libelous} exist, whence the proof follows.

\qed

  \textbf{Proof of  \cref{lem:slander}}. Conditionally on the trajectory $Y=\{Y_t, t\in[0,T]\}$, $X_T$ is a Gaussian random variable.
   Therefore, for any Borel set $A\subset \R$
  of zero Lebesgue measure, we have
  \[
   \mathrm{P} \{ X_T \in A \} =  \mathrm{E} (\I_{ X_T \in A })=
   \mathrm{ E}\left( \mathrm{E }\left( \I_{ X_T \in A }| \{Y_t, t\in[0,T]\} \right)\right)=0.
    \]

  The absolute continuity of the law of $S_T$
  follows from that of $X_T$ since
   $S _T = \exp\{X_T\}$.
   \qed

\textbf{Proof of \cref{lem:decry}}.
Let the function $H$ be locally Lipschitz
and $H'(x) = h(x)$ a.e.
with respect to the Lebesgue measure.
 Assume additionally that $h$ is of exponential growth. Hence it follows from   \cref{important rem} that

\[
 D ^W _u H(X_T) = h(X _T) D ^W _u X_T.
\]

Establish now  that $H(X_T) \in \mathrm{D}^{1,2}$, where we consider stochastic differentiation w.r.t. $W$, i.e., $\mathrm{D}^{1,2}=\mathrm{D}^{W,1,2}.$
Indeed,
  $h$ is of exponential  growth,
\[
 h(x)  \leq C_h (1+ e^{p_h|x|}),
\]
and
\[
 H(x)  = \int _0 ^x h(y)dy
 \leq C_h |x| (1+ e^{p_h|x|})\leq C_h   (1+ e^{(p_h+1)|x|}).
\]
Furthermore,
\[
 e^{(p_h+1)|X_T|} = (S_T) ^{p_h+1} \vee (S_T) ^{-p_h-1},
\]
we get from \cref{eq:exp-finity} that $\mathrm{E} H^2(X_T) < \infty$.
 Additionally,
\[
 \mathrm{E } \int\limits _0 ^T \left( h(X _T) D _u ^W X_T \right) ^2 du
 = \mathrm{E}
 \left( h ^2(X _T)
 \int\limits _0 ^T    \sigma ^2 (Y _u)   du\right)=
\]
\[
 \leq C  \left( \mathrm{E}h ^4(X _T)
 \int\limits _0 ^T \mathrm{E}\sigma ^4 (Y _u)   du\right) ^{1/2} < \infty.
\]
Therefore
\begin{equation}\label{eq:defame}
  H( X_T ) \in \text{D}^{1,2}.
\end{equation}

  Having established  both the existence and the form of the stochastic derivative, together with \cref{eq:defame}, we can proceed as in the proof of Proposition 4.1 \cite{AN15}.
Namely, the Skorokhod integral is the adjoint
operator to the Malliavin derivative, therefore
\begin{equation} \begin{gathered}\label{eq:din}
 \mathrm{E} h(X_T) = \frac 1T \mathrm{E }\left(
   \int\limits _0 ^T
   h(X_T ) D^W _u X_T \frac{1}{D^W_u X_T} du
 \right)
 =
 \frac 1T \mathrm{E} \left(
   \int\limits _0 ^T
   D^W _u H(X_T) \frac{1}{\sigma(Y_u)} du
 \right)\\
  = \frac 1T \mathrm{E}
 \left(
 H(X_T)
 \int\limits _0 ^T
   \frac{1}{\sigma(Y_u)} dW_u
 \right)=\frac 1T\mathrm{E}\left(
 H(X_T)Z(T)\right).
\end{gathered}\end{equation}
In particular, the function $G$ is  locally Lipschitz
and $G'(x) = g(x)$ a.e.
with respect to the Lebesgue measure. Moreover, $g$ is of exponential  growth, namely,
\[
 g(x)  \leq C_f (1+ e^{p|x|}),
\]
therefore   \cref{eq:for g} follows directly from \cref{eq:din}.

To establish \cref{eq:for f}, we start with the identity
\[
 G(x) = \frac{F(e^x)}{e^x} +
 \int\limits _0 ^x \frac{F(e^y)}{e^y}dy - F(1),
\]
then we rewrite it, applying  \cref{eq:for g}, as follows:
\begin{equation*} \begin{gathered}
 \mathrm{E} f(S_T) = \mathrm{E}g(X_T) = \frac 1T \mathrm{E}\left( G(X_T)Z_T \right)\\
= \frac 1T \mathrm{E}\left(  \left(\frac{F(S_T)}{S_T} +
 \int\limits _0 ^{X_T} \frac{F(e^y)}{e^y}dy - F(1)
 \right)  Z_T \right) = \frac 1T
 \mathrm{E}\left(\frac{F(S_T)}{S_T}   Z_T \right)\\+
 \frac 1T
 \mathrm{E}\left( \int\limits _0 ^{X_T} \frac{F(e^y)}{e^y}dy Z_T \right)-\frac 1T \mathrm{E}(F(1)Z_T)\\=\frac 1T
 \mathrm{E}\left(\frac{F(S_T)}{S_T}   Z_T \right)+
 \frac 1T
 \mathrm{E}\left( \int\limits _0 ^{X_T} \frac{F(e^y)}{e^y}dy Z_T \right).
\end{gathered}\end{equation*}
Applying equation \cref{eq:din} to $h(x)=  \frac{F(e^x)}{e^x}$, we get   that
\[
 \mathrm{E}\left(\frac{F(S_T)}{S_T} \right)= \frac 1T
 \mathrm{E}\left( \int\limits _0 ^{X_T} \frac{F(e^y)}{e^y}dy Z_T \right).
\]
Hence
\[
 \mathrm{E}f(S_T) = \frac 1T
 \mathrm{E}\left(\frac{F(S_T)}{S_T}   Z_T \right) +  \mathrm{E}\left(\frac{F(S_T)}{S_T}\right)
 = \mathrm{E}\left(
  \frac{F(S_T)}{S_T}
  \left(
  1 + \frac{Z_T}{T}
  \right)
  \right).
\]

\qed

 \textbf{Proof of   \cref{lem:disgruntle}}. $(i)$  Since the process $Y$ is Gaussian, it is sufficient to consider $\theta=2.$
 According to inequality  (1.9.2) from
 \cite{Mishura}, there exists a constant $C_H$ such that for any function $f\in L_{\frac1H}[0,T]$
 \begin{equation}\label{eq:1.9.2}
 \mathrm{E}\left(\int_0^T f(s)dB^H_s\right)^2\leq
 C_H\|f\|^2_{L_{\frac1H}[0,T]}.
 \end{equation}
Now, the increment of   $Y$ can be presented as
\[
Y_t - Y _s =  \int\limits _0 ^t
  h_{t,s}(\tau) d B^H(\tau) + Y_0 (e^{-\alpha t}
  - e^{-\alpha s}),
 \]
 where
 \[
  h_{t,s}(\tau) = \mathds{1} _{(s,t]} (u ) e^{- \alpha (t-u)}
  + \mathds{1} _{[0,s]} (u )  ( e^{- \alpha (t-u)}
  - e^{- \alpha (s-u)}).
 \]
 Note that
 \[
  |e^{- \alpha (t-u)}
  - e^{- \alpha (s-u)}| \leq \alpha (t-s), \ \
 | e^{- \alpha t}
  - e^{- \alpha s}| \leq \alpha (t-s),
 \]

and these simple inequalities imply, in particular, that
 \[
  \int\limits _0^t |h_{t,s}(u)|^{1/H} du =\int\limits _0^s|e^{- \alpha (t-u)}
  - e^{- \alpha (s-u)}|^{1/H} du+\int\limits _s^te^{- \frac{\alpha (t-u)}{H}}du \leq C |t-s|,
 \]
 and $(i)$   follows
 from \cref{eq:1.9.2}.\\
 $(ii)$ Again, since $Y$ and $Y^n$ both are Gaussian processes, it is sufficient to consider only $\theta=2.$
 Define the approximation
 $e_n (s) =
 e^{-\alpha   t_{i}}, \quad s \in [ t_i, t_{i+1} ),\; 0\leq i\leq n-1.$
  Then it follows from \cref{eq:1.9.2} that
 \[
  \mathrm{E}\left(Y_{t_j} -
    Y^n _{t_j} \right) ^{2}
    =
    \mathrm{E }\left[ \int\limits _0 ^{t_j} \Big(
    e^{-\alpha (t_j-s)} -
    e_n (t_j-s)
    \Big)dB^H _s
    \right]^{2}
 \]
 \[
  \leq C_H
  \left( \int\limits _0 ^{t_j} \Big(
    e^{-\alpha (t_j-s)} -
    e_n(t_j-s)
    \Big)^{1/H} ds
    \right)^{2H}
 \]
\[
 \leq
 C_H
  \left( \int\limits _0 ^{t_j} \alpha n^{-1/H} ds
    \right)^{2 H} =
    C_H  \alpha
     \left(
     t_j n^{-1/H}
     \right)^{2H} = Cn^{-2}.
\]
$(iii)$ Now, let
  $s \in [t_i, t_{i+1})$ and $\theta\geq 1$. Then it follows from $(i)$ and $(ii)$ that
\begin{equation*}\begin{gathered}
\mathrm{E}|Y _s - Y ^n_{s}|^{\theta} =
\mathrm{E}|Y _s - Y ^n_{t_i}|^{\theta}
\leq C \mathrm{E}|Y _s - Y _{t_i}|^{\theta}
+ C
\mathrm{E}|Y_{t_i} - Y ^n_{t_i}|^{\theta}\\
 \leq C \left(n^{-\theta H} +
 n^{-\theta}\right)\leq C n^{-\theta H}.
\end{gathered}\end{equation*}
Statement $(iv)$  follows  immediately from \cref{eq:1.9.2},
since
\[
 Y^n _{t_j} = \int\limits _{0} ^{t_j} e_n(t_j - s) d B^H _s,
\]
and functions $e_n$ are uniformly bounded in $n$.
\qed

   \textbf{Proof of   \cref{lem:tattle}}.
  Let us start with \cref{eq:haggle}.
  Taking into account condition $(\mathbf{B})$, $(ii)$ and $(iii)$, we can write
 \begin{equation}\begin{gathered} \label{eq:dazzle}
   \mathrm{E}(X_T - X ^n _T) ^2
   = \mathrm{E} \left[
  - \frac 12 \int\limits _0 ^T \sigma ^2 (Y _s)ds
 + \int\limits _0 ^T \sigma  (Y _s)dW _s
  +
 \frac 12 \int\limits _0 ^T \sigma ^2 (Y^n _s)  ds
 - \int\limits _0 ^T \sigma (Y^n _s)   dW _s
   \right] ^2\\
    \leq
  2   \mathrm{E} \left[
   \frac 12 \int\limits _0 ^T
   (\sigma ^2 (Y _s) - \sigma ^2 (Y^n _s) )ds
   \right]^2
   +
   2  \mathrm{E} \left[
   \int\limits _0 ^T (\sigma  (Y _s) -\sigma (Y^n _s)  ) dW _s
     \right]^2\\\leq
    \frac T2 \int\limits _0 ^T \mathrm{E}
   (\sigma ^2 (Y _s) - \sigma ^2 (Y^n _s) )^2ds
   +
   2
   \int\limits _0 ^T \mathrm{E} (\sigma  (Y _s) -\sigma (Y^n _s)  )^2 ds\\
    =\frac T2 \int\limits _0 ^T \mathrm{E}
   \left[ |\sigma  (Y _s) - \sigma (Y^n _s) |^2
   |\sigma  (Y _s) + \sigma (Y^n _s) |^2 \right]ds
   +2\int\limits _0 ^T \mathrm{E} (\sigma  (Y _s) -\sigma (Y^n _s)  )^2 ds
   \\
    \leq C
  \int\limits _{0} ^{T}  \Big(
  \mathrm{E} \left( \left| Y _s - Y ^n_{s} \right|^{2r}
  \left( C + |Y _s|^{2q} + |Y ^n_{s}|^{2q} \right)
  \right)\\
  + \mathrm{E} \left| Y _s - Y ^n_{s} \right|^{2r}
  \Big)ds
  \end{gathered}
\end{equation}
\[
 \leq C
  \int\limits _{0} ^{T}
  \left( \mathrm{E} \left| Y _s - Y ^n_{s} \right|^{4r}
 \mathrm{E} \left( C + |Y _s|^{4q } +
 |Y ^n_{s}|^{4q} \right)\right)^{1/2}
  ds.
\]
    \cref{lem:lemma 2.1} $(i)$, and    \cref{lem:disgruntle}  $(ii)$ imply that for any  $\theta\geq 1$
\begin{equation}\label{eq:sup in n}
\sup\limits _{n \in \N,
s \in [0,T]}\mathrm{E} \left(  |Y _s|^{\theta} +
|Y ^n_{s}|^{\theta} \right) < \infty.
\end{equation}
Moreover,
it follows from    \cref{lem:disgruntle} that
  for any  $s \in [0,T]$ and $\theta\geq 1$
\begin{equation}\begin{gathered}\label{eq:bound for diff-1}
\mathrm{E}|Y _s - Y ^n_{s}|^{\theta r}
 \leq C  n^{-\theta rH}.
\end{gathered}\end{equation}
Put $\theta=4q$ in \cref{eq:sup in n} and $\theta=4$ in \cref{eq:bound for diff-1}  and   substitute the result into  the right-hand side of \cref{eq:dazzle}:
\[
 \mathrm{E}(X_T - X ^n _T) ^2  \leq  C
 \int\limits _{0} ^{T}
  \left(\mathrm{E} \left( Y _s - Y ^n_s \right)^{4r}\right)^{\frac12}
 ds \leq C n^{-2rH},
\]
so that  \cref{eq:haggle} is proved.  Now continue with \cref{eq:bound for Z}.
Taking into account condition $(\mathbf{B})$, $(i)$, we get that  $$\left|\frac{1}{\sigma (x)} - \frac{1}{\sigma (y)}\right|
\leq \frac{|\sigma(x) - \sigma (y)|}{\sigma(x) \sigma (y)}
\leq \frac{|\sigma(x) - \sigma (y)|}{\sigma _{\min} ^2},$$
whence

\[
\mathrm{E}(Z_T - Z ^n _T) ^2
=
\int\limits _0 ^T
  \left(\frac{1}{\sigma (Y   )} -
  \frac{1}{\sigma (Y^n _s)   } \right)^2 ds
\]
\[
 \leq \frac {1}{\sigma _{\min} ^2} C_{\sigma}
  \int\limits _{0} ^{T}
  \mathrm{E} \left( Y _s -Y ^n_{s} \right)^{2r} ds.
\]
We can apply \cref{eq:bound for diff-1}   with $\theta=2$ to the latter inequality and conclude  this part of the proof   exactly as it was done  for \cref{eq:haggle}.

\qed

   \textbf{Proof of   \cref{lem:pull in horns}}.   We can write
  \begin{equation}\label{eq:vacillate}
    \mathrm{E} \left|\frac{F(S_T)}{S_T} - \frac{F(S _T ^n)}{S _T ^n}\right| ^2
    \leq
   2 \mathrm{E}\left|\frac{F(S_T)}{S_T} - \frac{F(S _T )}{S _T ^n}\right|^2
    +2 \mathrm{E}\left|\frac{F(S_T)}{S_T ^n} - \frac{F(S _T ^n)}{S _T ^n}\right|^2
    := 2I_1 + 2I_2.
  \end{equation}
 Now    we estimate  the right-hand side of \cref{eq:vacillate}  term by term. For $I_1$ we have that

\begin{equation} \label{eq:regurgitate}
  I_1 = \mathrm{E} \left( F(S_T) \big( ({S _T })^{-1}-({S _T ^n})^{-1} \big)\right)^2
  \leq
  \left( \mathrm{E} (F(S_T))^4  \mathrm{E} \big( ({S _T })^{-1}-({S _T ^n})^{-1} \big)^4
  \right) ^{1/2}.
\end{equation}

 On one hand, since $f$  consequently $F$ both have a polynomial growth, $\mathrm{E}( F(S_T))^4 < \infty$ according to   \cref{polyn growth}.  On the other hand,
 \begin{equation}\label{eq:gallivant}
   \mathrm{E} \big( ({S _T })^{-1}-({S _T ^n})^{-1} \big)^4 =
   S _0 ^{-4} e^{-4bT}
 \end{equation}
   \[
   \times \mathrm{E}
   \left(  \exp \left\{\frac12 \int\limits _0 ^{T} \sigma ^2 (Y^n _s)   ds
 - \int\limits _0 ^{T} \sigma (Y^n _s)    dW _s \right\}
 -
 \exp \left\{   \frac12  \int\limits _0 ^T \sigma ^2 (Y _s)ds
 - \int\limits _0 ^T \sigma  (Y _s)dW _s  \right\}
   \right)^4.
   \]

   Using the inequalities
   \[
    |e^x - e^y| \leq (e^x + e^y)|x-y|, \ \ \ x,y  \in \R,
   \]
    \[
    (x+y)^{2n} \leq C(n) (x ^{2n} + y ^{2n}), \ \ \ x,y  \in \R, n\in \N,
   \]
along with results outlined in    \cref{polyn growth} and   \cref{boundinn}, the Burkholder--Gundy  and H\"{o}lder inequalities, condition $(\mathbf{B})$, $(ii)$ and $(iii)$, and
  relation \cref{eq:sup in n} with $\nu=16q$, we get
   from \cref{eq:gallivant} that
   \begin{equation*}
\begin{gathered}
\mathrm{E} \big( ({S _T })^{-1}-({S _T ^n})^{-1} \big)^4\\
\leq
C\mathrm{ E}
\left( \exp \left\{\frac12 \int\limits _0 ^{T} \sigma ^2 (Y^n _s) ds
- \int\limits _0 ^{T} \sigma (Y^n _s) dW _s \right\}
-\exp \left\{ \frac12 \int\limits _0 ^T \sigma ^2 (Y _s)ds
- \int\limits _0 ^T \sigma (Y _s)dW _s \right\}
\right)^4\\
\leq C \mathrm{E} \left(
\left( \exp \left\{2 \int\limits _0 ^{T} \sigma ^2 (Y^n _s) ds
- 4 \int\limits _0 ^{T} \sigma (Y^n _s) dW _s \right\}
+
\exp \left\{ 2 \int\limits _0 ^T \sigma ^2 (Y _s)ds
- 4 \int\limits _0 ^T \sigma (Y _s)dW _s \right\}
\right) \right.
\\
\times \left.
\left( \frac12 \int\limits _0 ^{T} \sigma ^2 (Y^n _s) ds
- \int\limits _0 ^{T} \sigma (Y^n _s) dW _s
- \frac12 \int\limits _0 ^T \sigma ^2 (Y _s)ds
+ \int\limits _0 ^T \sigma (Y _s)dW _s
\right) ^4 \right)
\\
\leq C \left(\mathrm{ E}
\left( \exp \left\{4 \int\limits _0 ^{T} \sigma ^2 (Y^n _s) ds
- 8 \int\limits _0 ^{T} \sigma (Y^n _s) dW _s \right\}
\right. \right.
\\
+
\left. \left.
\exp \left\{ 4 \int\limits _0 ^T \sigma ^2 (Y _s)ds
- 8 \int\limits _0 ^T \sigma (Y _s)dW _s \right\}
\right) \right) ^{1/2}
\\
\times
\left[ \mathrm{E} \left( \frac12 \int\limits _0 ^{T} \sigma ^2 (Y^n _s) ds
- \int\limits _0 ^{T} \sigma (Y^n _s) dW _s
- \frac12 \int\limits _0 ^T \sigma ^2 (Y _s)ds
+ \int\limits _0 ^T \sigma (Y _s)dW _s
\right) ^8 \right]^{1/2}
\\
\leq C
\left[ \mathrm{E} \left( \frac12 \int\limits _0 ^{T} \sigma ^2 (Y^n _s) ds
- \int\limits _0 ^{T} \sigma (Y^n _s) dW _s
- \frac12 \int\limits _0 ^T \sigma ^2 (Y _s)ds
+\int\limits _0 ^T \sigma (Y _s)dW _s
\right) ^8 \right]^{1/2}
\\
\leq
C \left[ \mathrm{E} \left( \int\limits _0 ^{T} \sigma ^2 (Y^n _s) ds
- \int\limits _0 ^T \sigma ^2 (Y _s)ds \right) ^8
+
\mathrm{E} \left( \int\limits _0 ^{T} \sigma (Y^n _s) dW _s
- \int\limits _0 ^T \sigma (Y _s)dW _s \right)^8
\right]^{1/2}\\
\leq C
\left[
T ^7 \mathrm{E} \left( \int\limits _0 ^{T} \left( \sigma ^2 (Y^n _s)
- \sigma ^2 (Y _s)\right)^{8} ds \right)
+C \mathrm{ E}
\left( \int\limits _0 ^{T} \left( \sigma (Y^n _s)
- \sigma (Y _s)\right) ^{2} ds \right)^4
\right]^{1/2}\\
= C
\left[
T ^7 \left( \int\limits _0 ^{T} \mathrm{E} \left\{ (\sigma (Y^n _s)
- \sigma (Y _s))
(\sigma (Y^n _s) +
\sigma (Y _s))
\right\}^{8} ds \right)
\right. \\
\left.
+
CT^3 \mathrm{E}
\left( \int\limits _0 ^{T} \left( \sigma (Y^n _s)
- \sigma (Y _s)\right) ^{8} ds \right)
\right]^{1/2}\end{gathered}\end{equation*}
\begin{equation}\begin{gathered}\label{eq:go out of your way}
\leq
    C
   \left[
  \mathrm{E }\left(   \int\limits _0 ^{T} \left(
  \mathrm{E}|Y _s - Y ^n _s|^{16r}
 \mathrm{E}\left(C + |Y _s|^{16q} + |Y ^n _s|^{16q}\right)
 \right)^{1/2} ds \right) \right.\\
 + \left. \mathrm{E}
 \left(  \int\limits _0 ^{T}
 \left| Y _s - Y ^n _s \right|^{8r}  ds \right)
    \right]^{1/2}\\ \leq C
 \left[
   \int\limits _0 ^{T}  \left(  \left(\mathrm{E}
   \left| Y _s - Y ^n _s \right|^{16r}\right)^{\frac{1}{2}}
   + \mathrm{E }\left| Y _s - Y ^n _s \right|^{8r} \right) ds
 \right]^{1/2}.
\end{gathered}\end{equation}

Applying   \cref{eq:bound for diff-1} consequently with $\theta=8$ and $\theta=16$
we get that the last expression in \cref{eq:go out of your way} does not exceed
$
 C \left( \frac 1n \right) ^{4rH},
$
thus from \cref{eq:regurgitate} we obtain
\begin{equation}\label{eq:I1}
 I_1 \leq Cn^{-2rH}.
\end{equation}

Now we continue with $I_2$ from the relation \cref{eq:vacillate}:
 \[
  I_2 \leq \left[
  \mathrm{E}( {F(S_T)} - {F(S _T ^n)} )^4
  \right] ^{1/2} \left[
  \mathrm{E} ({S _T ^n})^{-4}
  \right] ^{1/2}.
\]
The second multiplier is bounded according to   \cref{boundinn}, therefore it follows from condition $(A)$, $(i)$, that
\begin{equation*}\begin{gathered}
 I_2 \leq C \left[
 \mathrm{E}( {F(S_T)} - {F(S _T ^n)} )^4
 \right] ^{1/2}
 =
 C \left[ \mathrm{E} \left( \int\limits _{S_T \wedge S _T ^n} ^{S_T \vee S _T ^n}
 f(x) dx \right)^4 \right] ^{1/2}\\
 \leq C(C_f )^2\left[ \mathrm{E}\left( |S_T - S _T ^n|^4(1+
 S_T ^p + (S _T ^n) ^p)^4 \right) \right] ^{1/2}
 \leq C \left[ \mathrm{E} |S_T - S _T ^n|^8
 \mathrm{E}(1+
 S_T ^p + (S _T ^n) ^p)^8
 \right] ^{1/4}.
\end{gathered}\end{equation*}

 According to   \cref{lem:lemma 2.1} and   \cref{boundinn},
\[\sup\limits _{n \in \N} \mathrm{E}(1+
 S_T ^p + (S _T ^n) ^p)^8 < \infty, \]
whence we get that

\[
 I _2 \leq  C \left[ \mathrm{E} |S_T - S _T ^n|^8
 \right] ^{1/4}.
\]
To evaluate the right-hand side of this inequality, we can proceed as in the proof
  of  \cref{eq:go out of your way} and subsequent inequalities,
because neither the opposite sign of the exponents
nor the 8th power instead of the 4th
  lead to serious discrepancies
in   the estimations. Therefore we get
\begin{equation}\label{eq:I2}
 I _2 \leq C
 \left[
   \int\limits _0 ^{T}  \left(  \left(\mathrm{E}
   \left| Y _s - Y ^n _s \right|^{32r}\right)^{\frac{1}{2}}
   + \mathrm{E }\left| Y _s - Y ^n _s \right|^{16r} \right) ds
 \right]^{\frac{1}{8}}  \leq  C  \left(
 n^{-16rH} \right)^{1/8} =
  C n^{-2rH}.
\end{equation}
Bounds \cref{eq:I1} and \cref{eq:I2} complete the proof.
  \qed

   \textbf{Proof of \cref{thm:convergence thm}}.
  By
     \cref{lem:decry}
   we can write
   \[
   \left|\mathrm{E} f(X_T) - \mathrm{E} \left(
  \frac{F(S ^n _T)}{S ^n _T}
  \left(
  1 + \frac {Z ^n _T}{T}
  \right)
  \right) \right|
  =
    \mathrm{E} \left| \left(
  \frac{F(S_T)}{S_T}
  \left(
  1 + \frac {Z  _T}{T}
  \right)
  \right)
  - \left(
  \frac{F(S_T^n)}{S_T^n}
  \left(
  1 + \frac {Z ^n _T}{T}
  \right)
  \right)
  \right|
 \]
\[
 \leq \frac 1T  \mathrm{E}\left|
  \frac{F(S_T)}{S_T}
  \left(
  Z_T - Z^n_T
  \right)
  \right|
  + \mathrm{E} \left| \left(
  1 + \frac {Z ^n _T}{T}
  \right)\left(
  \frac{F(S_T)}{S_T} -   \frac{F(S_T^n)}{S_T^n} \right) \right|
\]
 \[
   \leq \frac 1T
  \left[
  \mathrm{E} \left(
  \frac{F(S_T)}{S_T}   \right)^2
  \mathrm{E} \left(
  Z_T - Z^n_T
  \right)^2
  \right] ^{1/2}
 + \left[ \mathrm{E} \left(
  \frac{F(S_T)}{S_T} -  \frac{F(S_T^n)}{S_T^n}   \right)^2
  \mathrm{E}\left(
  1 + \frac {Z ^n _T}{T}
  \right)^2  \right] ^{1/2}.
   \]

  According to   \cref{lem:lemma 2.1},  \cref{polyn growth} and Cauchy--Schwartz inequality,   $\mathrm{E} \left(
  \frac{F(S_T)}{S_T}   \right)^2 < \infty$.
  Obviously,  $\sup_{n\geq 1}\mathrm{E}\left(
     {Z ^n _T}
  \right)^2 < \frac{T}{ \sigma_{min}^2}.$
  Now the proof follows from   \cref{lem:tattle} and
   \cref{lem:pull in horns}.
  \qed

 \textbf{Proof of  \cref{lem:lem 5.1}}. Proof immediately follows from the general formula for the density of $k$-dimensional Gaussian vector:
 \begin{equation}\label{eq:density}p(\bar{x})=(2\pi)^{-\frac{k}{2}}|C|^{-1}\exp\{-(C^{-1}(\bar{x}-\bar{a}), \bar{x}-\bar{a}\},\end{equation}
 where $\bar{x}\in \R^k$, $\bar{a}$ is a vector of expectations, $C$ is a covariance matrix. In our case covariance matrix equals \begin{equation*}
    C=C_{X,Z} = \begin{pmatrix}
    \sigma_{_Y} ^2 & T\\
    T & \sigma_{_Z} ^2
  \end{pmatrix},
   \end{equation*} $$k=2,\; |C_{X,Z}|=\Delta=\sigma_{Y} ^2\sigma _{Z} ^2-T^2, \bar{a}=(m_X,0)=(\log S_0+bT-\frac{1}{2}\sigma_{Y} ^2, 0),$$   and  \cref{eq:boink} follows immediately from \cref{eq:density}.

 \qed

 \textbf{Proof of  \cref{thm:analytic expression thm}}.  Applying   conditioning on $Y$, \cref{eq:for g}, and   \cref{lem:lem 5.1}, we get  that

\begin{equation}\begin{gathered}\label{eq:errand}
 T \mathrm{E} g(X_T) =
   \mathrm{E}\left( G(X_T)
 \int\limits _0 ^T
   \frac{1}{\sigma(Y_u)} dW_u
 \right)
 =
 \mathrm{E} \left(\mathrm{E}
 \left(
 G(X_T)
 \int\limits _0 ^T
   \frac{1}{\sigma(Y_u)} dW_u
   \bigg| \{Y_s, s \in [0,T]\}
 \right)\right)\\
 = \mathrm{E} \left(\mathrm{E}
 \left(
\int\limits _{\R ^2} G(x) z p_{X,Z}(x,z) dx dz
   \bigg| \{Y_s, s \in [0,T]\}
 \right)\right)=
 \mathrm{E} \int\limits _{\R ^2} G(x) z p_{X,Z}(x,z) dx dz\\
   =
  \mathrm{E} \int\limits _{\R } G(x)\left(\int\limits _{\R } z p_{X,Z}(x,z) dz\right) dx.
\end{gathered}\end{equation}

 The inner integral can be significantly  simplified. Indeed,
denote   $\tilde{x}= x-m_{_Y}$. Then
 \[
  \int\limits _{\R } z p_{X,Z}(x,z) dz
  = \frac{1}{2 \pi \Delta^{\frac{1}{2}}}
    \int\limits _{\R } z
    \exp
    \Bigg\{
   - \frac{1 }{2
   \Delta
   }
   \left(
    {\sigma _{_Z} ^2}{\tilde x ^2}
    +{\sigma _{_Y} ^2}{z^2}
    -
    {2T \tilde x z }
   \right)
    \Bigg\}
    dz
\]
\[
 =\frac{1}{2 \pi \Delta^{\frac{1}{2}}}
    \int\limits _{\R } z
    \exp
    \Bigg\{
   - \frac{1 }{2
   \Delta
   }
   \left(
  \left( \sigma _{_Y} z - \frac{T \tilde x}{  \sigma _{_Y}}  \right) ^2
   - \frac{T^2 \tilde x ^2 }{  \sigma _{_Y} ^2}
   + \sigma _{_Z} ^2 \tilde x ^2
   \right)
    \Bigg\}
    dz
\]
\[
 = \frac{1}{2 \pi \Delta^{\frac{1}{2}}}
 \exp
    \Bigg\{
    -\frac{\tilde x ^2}{2\Delta}
     \frac{\sigma _{_Y} ^2 \sigma _{_Z} ^2
 - T ^2}{\sigma _{_Y} ^2}
       \Bigg\}
    \int\limits _{\R } z
    \exp \Bigg\{
   - \frac{1 }{2
   \Delta
   }
   \left( \sigma _{_Y} z - \frac{T \tilde x}{  \sigma _{_Y}} \right) ^2
   \Bigg\}
    dz
\]
\[
 =
 \frac{1}{2 \pi \Delta^{\frac{1}{2}}}
 \exp
    \Bigg\{
    -\frac{\tilde x ^2}{2\sigma _{_Y} ^2}
       \Bigg\}
    \int\limits _{\R } z
 \exp \Bigg\{
   -
  \left( \frac{\sigma _{_Y}}{\sqrt{2}\Delta^{\frac{1}{2}}} z
   - \frac{T \tilde x}{ \sqrt{2} \Delta^{\frac{1}{2}} \sigma _{_Y}}  \right) ^2
   \Bigg\}
    dz.
\]

Since
\[
 \int\limits _{\R } x e^{-(ax - b)^2} dx = \frac{b}{a^2}\sqrt{\pi},
\]
we obtain

\begin{equation}\label{eq:blow your wad}
  \int\limits _{\R } z p_{X,Z}(x,z) dz
       =
       \frac{T\tilde x}{\sigma^3_{Y} \sqrt{2\pi}}
        \exp  \Bigg\{
    -\frac{\tilde x ^2}{2\sigma _{_Y} ^2}
       \Bigg\}.
    \end{equation}

Combining \cref{eq:errand}
and \cref{eq:blow your wad},
we get   the proof.

\qed

\textbf{Proof of   \cref{thm:second conv thm}}. To simplify notations, without loss of generality, let us assume that  $X_0+bT=0.$
Then, using \cref{eq:assuage}, we get  that

\begin{equation*}
 \begin{gathered}
  \Bigg| \E g(X _T) - (2\pi)^{ - \frac{1}{2}}\int\limits _{\R } G(x)\mathrm{E}\Bigg(
 \frac { (x - m _{_{Y,n}}) }{\sigma_{_{Y,n}} ^3}
        \exp  \Bigg\{
    -\frac{ (x - m _{_{Y,n}}) ^2}{2\sigma _{_{Y,n}} ^2}
       \Bigg\} \Bigg)dx \Bigg|
\\
 =(2\pi)^{- \frac{1}{2}}\Bigg|\int\limits _{\R } G(x)
 \mathrm{E}\Bigg(
 \frac { (x - m _{_Y}) }{\sigma_{_Y}^3}
        \exp  \Bigg\{
    -\frac{ (x - m _{_Y}) ^2}{2\sigma _{_Y} ^2}
       \Bigg\} -
      \\ \frac { (x - m _{_{Y,n}}) }{\sigma_{_{Y,n}} ^3}
        \exp  \Bigg\{
    -\frac{ (x - m _{_{Y,n}}) ^2}{2\sigma _{_{Y,n}} ^2}
       \Bigg\}
       \Bigg)
       dx
       \Bigg|
\\ 
 \leq
 (2\pi)^{ - \frac{1}{2}}\Bigg(\int\limits _{\R } G(x)\Bigg(
 \mathrm{E}\Bigg(\Bigg|
 \frac { x - m _{_Y} }{\sigma_{_Y}^3}
 -
 \frac {  x - m _{_{Y,n}}  }{\sigma_{_{Y,n}} ^3}
  \Bigg|
        \exp  \Bigg\{
    -\frac{ (x - m _{_Y}) ^2}{2\sigma _{_Y} ^2}
       \Bigg\}\Bigg)
       \\  +
  \mathrm{E}\Bigg|
 \frac { (x - m _{_{Y,n}}) }{\sigma_{_{Y,n}} ^3}
 \Bigg(
 \exp  \Bigg\{
    -\frac{ (x - m _{_Y}) ^2}{2\sigma _{_Y} ^2}
       \Bigg\}
 -  \exp  \Bigg\{
    -\frac{ (x - m _{_{Y,n}}) ^2}{2\sigma _{_{Y,n}} ^2}
       \Bigg\}
\Bigg)
       \Bigg|
       \Bigg)dx
       \Bigg)
\\
 :=  (2\pi)^{ - \frac{1}{2}}\Bigg(\int\limits _{\R } G(x)(J_1(x) + J_2(x))dx\Bigg).
 \end{gathered}
\end{equation*}

To estimate $J_1(x)$, denote $E_{exp}(x)=\left(\E
 \exp  \Bigg\{
    - \frac{ (x - m _{_Y}) ^2}{\sigma _{_Y} ^2}
       \Bigg\}
        \right)^{1/2}$    and notice that

\[
  \frac {  x - m _{_Y}  }{\sigma_{_Y} ^3}
 -
 \frac { x - m _{_{Y,n}}  }{\sigma_{_{Y,n}} ^3}
 = \sigma_{_Y} ^{-3}(m _{_{Y,n}} - m _{_Y})
 + (x - m _{_{Y,n}})(\sigma_{_Y}^{-3}-\sigma_{_{Y,n}} ^{-3}).
\]
Hence
\begin{equation} \label{eq:morass}
 \begin{gathered}
   J_1(x) \leq C \left(
 \E\left(
 \sigma_{_Y} ^{-3}(m _{_{Y,n}} - m _{_Y})
 \right)^2
 +
 \E\left(
 (x - m _{_{Y,n}})(\sigma_{_Y}  ^{-3}-\sigma_{_{Y,n}} ^{-3})
 \right)^2 \right)^{1/2}E_{exp}(x).
  \end{gathered}
\end{equation}

Since $|a_1^3 - a_2 ^3|
\leq |a_1^2 - a_2^2|(a_1+a_2)$, $a_1, a_2 >0$,
and also   the lower bounds
$ \sigma_{_{Y}}^{2} \geq T \sigma _{\min}^2 $,
$\sigma_{_{Y,n}}^{2} \geq T \sigma _{\min}^2 $ hold,
one  can conclude that
\begin{equation}\label{eq:harrowing}
\begin{gathered}
 \left| \sigma_{_Y}  ^{-3}-\sigma_{_{Y,n}} ^{-3} \right|
 \leq \left| \sigma_{_Y}  ^{-2}-\sigma_{_{Y,n}} ^{-2} \right|
 \left( \sigma_{_{Y,n}} ^{-1}+\sigma_{_Y} ^{-1}\right)
 =
 \frac{\left| \sigma_{_Y} ^{2}-\sigma_{_{Y,n}} ^{2} \right|}
 {\sigma_{_Y} ^{2}\sigma_{_{Y,n}} ^{2}}\left( \sigma_{_{Y,n}} ^{-1}
 +\sigma_{_Y} ^{-1}\right)
 \\
 \leq \frac{2\left| \sigma_{_Y} ^2-\sigma_{_{Y,n}} ^{2} \right|}
 {\sigma_{_{Y,n}}^2 T^{\frac32} \sigma _{\min}^3}.
\end{gathered}
\end{equation}
Therefore \begin{equation*}
\begin{gathered}|(x - m _{_{Y,n}})(\sigma_{_Y}  ^{-3}-\sigma_{_{Y,n}} ^{-3})|\leq \frac{2|x +\frac12 \sigma_{_{Y,n}} ^{2}| | \sigma_{_Y}  ^{2}-\sigma_{_{Y,n}} ^{2} |}{\sigma_{_{Y,n}^2} T^{\frac32} \sigma_{\min}^3}\leq C(1+|x|)\left| \sigma_{Y} ^{2}-\sigma_{Y,n} ^{2} \right|.
\end{gathered}
\end{equation*}

Since $m _{_{Y,n}} - m _{_Y} = -
\frac 12 (\sigma_{_{Y,n}}^2 -\sigma_{_Y} ^2)$,
we get from  \cref{eq:morass} and \cref{eq:harrowing}
that

       \begin{equation}\label{eq:goon}
        \begin{gathered}
         J_1(x)
 \leq C(1+|x|)( \E (\sigma_{_{Y,n}} ^2 -\sigma_{_Y} ^2 )^2) ^{1/2}E_{exp}(x).
  \end{gathered}
       \end{equation}

Similarly to \cref{eq:dazzle} and \cref{eq:go out of your way}, we get, applying condition $(\mathbf{B})$,   \cref{lem:disgruntle}, $(iii)$ and $(iv)$, together with  the standard  H\"older's inequality, that
\begin{equation}\begin{gathered}\label{eq:bound-sigma}
 \E (\sigma_{_{Y,n}}^2 - \sigma_{_Y}^2)^2
  = \E \left(\int\limits _0 ^T
    (\sigma^2 (Y ^n _s) -
    \sigma^2 (Y  _s))ds
    \right)^2 \leq
 T \E \int\limits _0 ^T(
    \sigma^2 (Y ^n _s) -
    \sigma^2 (Y  _s) )^2ds
\\
 \leq C_{\sigma}C \int\limits _0 ^T
 \left[
     \E (Y ^n _s -   Y  _s)^{4r}
   \E  (\sigma (Y ^n _s) + \sigma (Y  _s))^{4}
 \right]^{1/2}ds
 \leq
 C \int\limits _0 ^T
 \left[
     \E (Y ^n _s -   Y  _s)^{4r}
 \right]^{1/2}ds \leq C n^{-2rH}.
\end{gathered}\end{equation}

Combining the latter inequality  with \cref{eq:goon} we get that

\begin{equation} \label{eq:tawdry}
   J_1(x)  \leq
 C  n^{-rH} (1+|x|)
 E_{exp}(x),
\end{equation}
and consequently
\begin{equation}\label{eq:ire}
\int\limits _{\R } G(x) J_1(x)  dx\leq
C n^{-rH}
     \int\limits _{\R } G(x) (1+|x|)
E_{exp}(x)
       dx.
\end{equation}
Let us show that  the integral in the right--hand side of \cref{eq:ire} is finite. Applying the standard H\"{o}lder inequality together with polynomial growth of $G(x)$, we get that
\begin{equation}\begin{gathered}\label{eq:J-1}
 \int\limits _{\R } G(x) (1+|x|)
 E_{exp}(x)
       dx\leq \left(\int\limits _{\R } G^2(x) (1+|x|)^2e^{-(2p+1)|x|}dx\right)^{\frac{1}{2}}
       \\
       \times\left(\int\limits _{\R }e^{(2p+1)|x|}\E
 \exp  \Bigg\{
    - \frac{ (x - m _{_Y}) ^2}{\sigma _{_Y} ^2}
       \Bigg\}dx\right)^{1/2}
       \\
       \leq  C \left( 2\int\limits _0 ^\infty e^{(2p+2)x}\E
 \exp  \Bigg\{
    - \frac{ x ^2}{\sigma _{_Y} ^2} - \frac{\sigma _{_Y} ^2}{4}
       \Bigg\}dx\right)^{1/2}
       \\
       =
       C \left( 2 \E\int\limits _0 ^\infty
 \exp  \Bigg\{
    - \left(  \frac{x}{\sigma _{_Y}} - (p+1)\sigma _{_Y} \right)^2
    + (p+1)^2\sigma _{_Y} ^2 - \frac{\sigma _{_Y} ^2}{4}
       \Bigg\}dx\right)^{1/2}
       \\
       =
      C \left( 2 \E \left[ \exp  \Bigg\{  (p+1)^2\sigma _{_Y} ^2 - \frac{\sigma _{_Y} ^2}{4}
      \Bigg\} \int\limits _0 ^\infty
 \exp  \Bigg\{
    - \left(  \frac{x}{\sigma _{_Y}} - (p+1)\sigma _{_Y} \right)^2
       \Bigg\}dx \right] \right)^{1/2}
       \\
       =
        C \left( 2 \E \left[ \exp  \Bigg\{  (p+1)^2\sigma _{_Y} ^2 - \frac{\sigma _{_Y} ^2}{4}
      \Bigg\} \int\limits _0 ^\infty
 \exp  \Bigg\{
    - \left(  \frac{x}{\sigma _{_Y}}  \right)^2
       \Bigg\}dx \right] \right)^{1/2}
       \\
       \leq
       C \left( \E \left[ \sigma _{_Y}
       \exp  \Bigg\{  (p+1)^2\sigma _{_Y} ^2 - \frac{\sigma _{_Y} ^2}{4}
      \Bigg\}
       \right]
       \right)^{1/2}.
      \end{gathered}\end{equation}
The finiteness of the expectation on the last line
follows from  condition $(\mathbf{B})$ and   \cref{lem:lemma 2.1},
 formula \cref{eq:exp-mom-sup},  since
 \[
 \sigma _{_Y} ^2 \leq C (1 + \sup\limits _{t \in [0,T]} |Y_t|^{2q}).
 \]
Construction of upper bound for  $J_2(x)$ is similar.
 Indeed, $$|\exp\{-u^2\}-\exp\{-v^2\}|\leq 2(|u|+|v|)(\exp\{-u^2\}+\exp\{-v^2\})|u-v|.$$ In our case
  $|u|=\left|\frac{x - m _{_{Y}}}{\sigma_{_{Y}}}\right|\leq C(1+|x|)\sigma_{_{Y}}$, $|v|=\left|\frac{x - m _{_{Y_n}}}{\sigma_{_{Y_n}}}\right|\leq C(1+|x|)\sigma_{_{Y_n}}$, and $\left|
 \frac { (x - m _{_{Y,n}}) }{\sigma_{_{Y,n}} ^3}\right|\leq C(1+|x|)$, therefore

\begin{equation}\label{eq:pylon}
 \begin{gathered}
  J_2 (x)\leq C(1+|x|)^2\mathrm{E}\Bigg((\sigma _{_{Y,n}}+\sigma _{_{Y}})\Bigg(\exp  \Bigg\{
    -\frac{ (x - m _{_Y}) ^2}{2\sigma _{_Y} ^2}\Bigg\}+
 \exp  \Bigg\{
    -\frac{ (x - m _{_{Y,n}}) ^2}{2\sigma _{_{Y,n}} ^2}
       \Bigg\}
 \Bigg)
 \\
   \times
  \left|\frac{  x - m _{_Y} }{ \sigma _{_Y}  }
 -\frac{  x - m _{_{Y,n}} }{ \sigma _{_{Y,n}} }
   \right|\Bigg)
  \leq  C(1+|x|)^2\left(\mathrm{E}\left(\frac{  x - m _{_Y}   }{ \sigma _{_Y}  }
 -
    \frac{  x - m _{_{Y,n}} }{ \sigma _{_{Y,n}} }\right)^2\right)^{\frac12}\\ \times\left(\mathrm{E}\left((\sigma _{_Y}+\sigma _{_Y,n})\left(\exp  \Bigg\{
    -\frac{ (x - m _{_Y}) ^2}{2\sigma _{_Y} ^2}
       \Bigg\}+
 \exp  \Bigg\{
    -\frac{ (x - m _{_{Y,n}}) ^2}{2\sigma _{_{Y,n}} ^2}
       \Bigg\}\right)\right)^2\right)^{\frac12}.
 \end{gathered}
\end{equation}
Now,  taking into account that $\sigma_{_Y}\wedge\sigma _{_{Y,n}}\geq T^{\frac12}\sigma_{min}$, and $\sigma _{_Y} \sigma _{_Y,n}(\sigma _{_Y}+\sigma _{_Y,n})\geq C\sigma^2 _{_Y,n}$, we get
\begin{equation*}\begin{gathered}\left|\frac{  x - m _{_Y}   }{ \sigma _{_Y}  } -\frac{  x - m _{_{Y,n}} }{ \sigma _{_{Y,n}} }\right|\leq C|\sigma _{_{Y,n}}^2-\sigma _{_{Y,n}}^2|+\frac{| x - m _{_{Y,n}}|}{ \sigma _{_Y} \sigma _{_Y,n}(\sigma _{_Y}+\sigma _{_Y,n})}|\sigma _{_Y}^2-\sigma _{_Y,n}^2|\\ \leq C(1+|x|)|\sigma _{_Y}^2-\sigma _{_Y,n}^2|\leq C(1+|x|)|\sigma _{_Y}^2-\sigma _{_Y,n}^2|,\end{gathered}\end{equation*}
and from \cref{eq:bound-sigma} we deduce that
$\left(\mathrm{E}\left(\frac{  x - m _{_Y}   }{ \sigma _{_Y}  }
 -\frac{  x - m _{_{Y,n}} }{ \sigma _{_{Y,n}} }\right)^2\right)^{\frac12}\leq C(1+|x|) n^{-rH}$.
 Together with \cref{eq:pylon} this implies that
     \begin{equation*}\begin{gathered}
     \int_{\R}  G(x) J_2 (x) dx \leq
    C n^{-rH} \int_{\R} (1+|x|)^3 G(x)\Bigg(\mathrm{E}\Bigg((\sigma _{_Y}+\sigma _{_Y,n})\\
     \times\left(\exp  \Bigg\{
    -\frac{ (x - m _{_Y}) ^2}{2\sigma _{_Y} ^2}
       \Bigg\}+
 \exp  \Bigg\{
    -\frac{ (x - m _{_{Y,n}}) ^2}{2\sigma _{_{Y,n}} ^2}
       \Bigg\}\right)\Bigg)^2\Bigg)^{\frac12}dx.
       \end{gathered}\end{equation*}
     The fact that  the   integral in the right-hand side    is finite,  can be established   via   the  same approach  as applied to the integral $\int\limits _{\R } G(x) (1+|x|)E_{exp}(x)dx$ in \cref{eq:J-1}.

\qed

 \textbf{Proof of  \cref{lem:negative moments}}. For any $\delta>0$ choose $\varepsilon=\varepsilon(\delta)$ in such a way that $\varepsilon(2+\alpha|Y_0|)<\delta$. Then we get from the representation \cref{eq:OU direct} that for $0\leq t\leq \tau_\varepsilon\wedge\varepsilon$
 \begin{equation*}\begin{gathered}|{Y_t} - {Y_0}|\leq |Y_0|(1-e^{-\alpha t})+|B_t^H|+\alpha e^{-\alpha t}\int_0^te^{\alpha s}|B_s^H|ds\leq  |Y_0|\alpha \varepsilon+\varepsilon+\varepsilon e^{-\alpha t}(e^{\alpha t}-1)\\ \leq \varepsilon(2+\alpha|Y_0|)<\delta.\end{gathered}\end{equation*}

Therefore for $\varepsilon<\frac{\delta}{2+\alpha|Y_0|}$ we have that $\nu_\delta>\tau_\varepsilon\wedge\varepsilon$. So, it is sufficient to prove that
for any $\varepsilon>0$ and any  $l>0$ $$\E(\tau_\varepsilon\wedge\varepsilon)^{-l}<\infty.$$
Now, for $v<\varepsilon$ $$\mathrm{P}\{\tau_\varepsilon\wedge\varepsilon<v\}=\mathrm{P}\{\tau_\varepsilon<v\}=\mathrm{P}\{\sup_{0\leq t\leq v}|B_t^H|\geq \varepsilon\}.$$
Furthermore, it follows from self-similarity and symmetry of the fBm that $\mathrm{P}\{\sup_{0\leq t\leq v}|B_t^H|\geq \varepsilon\}
\leq 2\mathrm{P}\{\sup_{0\leq t\leq 1} B_t^H \geq \frac{\varepsilon}{v^H}\}.$ Moreover,   denote $\vartheta=\E\sup_{0\leq t\leq 1} B_t^H$. Then, according to inequality (2.2) from \cite{Tal94} that for $\frac{\varepsilon}{v^H}>\vartheta$
$$\mathrm{P}\{\sup_{0\leq t\leq v} B_t^H \geq \frac{\varepsilon}{v^H}\}\leq \exp\left\{-\frac{\left(\frac{\varepsilon}{v^H}-\vartheta\right)^2}{2}\right\}=
\exp\left\{-\frac{( \varepsilon -\vartheta v^H)^2}{2v^{2H}}\right\},$$
whence the proof immediately follows.

 \qed
\begin{remark} \label{negative moments-1}  Exponential bounds for the distribution of $\tau_\varepsilon$ allow to prove that $\E(\tau_\varepsilon\wedge\varepsilon\wedge a)^{-l}<\infty$ for any $a,l>0.$\end{remark}

 \textbf{Proof of  \cref{lem:lurid}}.  As it follows from  Proposition 2.1.1 and Exercise 2.1.1 in \cite{Nua06}, it is sufficient to show that
  \begin{equation}\label{eq:repugnant}
   \sigma ^2 _{_Y} \in \mathrm{D} ^{2,4}
  \end{equation}
  and that
  \begin{equation}\label{eq:flip out}
   \mathrm{E} \left( ||D^B \sigma ^2 _{_Y}||_{H} \right) ^{-8} < \infty.
  \end{equation}
Recall that $\kappa(x)=\sigma(x)\sigma'(x)$. It follows from conditions $(\mathbf{B})$ and $(\mathbf{D})$ that $\kappa$ and $\kappa'$ are functions of polynomial growth, $\kappa(x)>0.$ Recall the notation  $l(u,s)=c _He^{-\alpha s}
  \int\limits _u ^ s e^{\alpha v} v ^{H - 1/2} (v-u) ^{H - 3/2}dv$. Taking into account \cref{eq:derivative} and \cref{eq:vertigo}, we write    the stochastic derivative as
\begin{equation*}
\begin{gathered}
  D^B_u ( \sigma ^2 _{_Y}) =
 D^B _u(   \int\limits _0 ^T
\sigma ^2(Y _s)ds) =  2
\int\limits _0 ^T \kappa(Y _s)
D^B _u  Y _s ds\\=  2c _H u ^{1/2 - H}
\int\limits _u ^T \kappa(Y _s)
 e^{-\alpha s}
  \int\limits _u ^ s e^{\alpha v} v ^{H - 1/2} (v-u) ^{H - 3/2}dvds\\=2 u ^{1/2 - H}\int\limits _u ^T \kappa(Y _s)
 l(u,s)ds.
\end{gathered}\end{equation*}
 Therefore, the iterated derivative equals
 \begin{equation}\label{eq:repeated}\begin{gathered}
  D^B_z(D^B_u ( \sigma ^2 _{_Y}))
   \\=  2  u ^{1/2 - H}z^{1/2 - H}
\int\limits _{u\vee z} ^T \kappa'(Y _s) l(z,s)
  l(u,s)ds.
\end{gathered}\end{equation}
  Obviously, the right-hand side of \cref{eq:repeated}   is in $H\otimes H$, and the corresponding integral has moments of any order, due to polynomial growth of $\kappa'$, which implies  \cref{eq:repugnant}.

  To prove \cref{eq:flip out},   note that

$$ D _u ^B (\sigma^2 _{_Y}) \geq C \int\limits _u ^T
 \kappa (Y_s)   (s-u) ^{H - 1/2} ds,
$$
whence
\[
 ||D^B  \sigma ^2 _{_Y}|| _{H} ^2 = \int\limits
 _{0} ^T \left( D _u ^B \sigma ^2 _{_Y} \right)^2 du
 \geq C \int\limits _{0} ^{ T} du
 \left( \int\limits _{u} ^{ T}  \kappa (Y_s) (s-u) ^{H - 1/2} ds \right)^2.
\]
Now, let $\sigma'(Y_0)=\sigma_0>0$. Choose $\delta>0$ so that for  $y\in[Y_0-\delta, Y_0+\delta]$ to provide on this interval the lower bound $\sigma'(y)>\frac{\sigma_0}{2}$. Then choose $\varepsilon=\varepsilon(\delta)$, as it was mentioned in the proof of   \cref{lem:negative moments}, and put $\zeta=\tau _{\varepsilon} \wedge \varepsilon\wedge\frac{T}{2}.$

Then
\[\int\limits _{0} ^{ T} du
 \left( \int\limits _{u} ^{ T}  \kappa (Y_s) (s-u) ^{H - 1/2} ds \right)^2
 \geq C \int\limits _{0}
 ^{\frac 13 \zeta}
 du
 \left( \int\limits
 _{ \frac 23 \zeta}
 ^{\zeta}
  \kappa (Y_s) (s-u) ^{H - 1/2} ds \right)^2
\]
\[
  \geq
 C \int\limits _{0}
 ^{\frac 13 \zeta}
 du
 \left( \int\limits
 _{ \frac 23 \zeta}
 ^{\zeta}
 \sigma  _{\min} \sigma_0
 \left(\frac 13 \zeta
 \right) ^{H - 1/2} ds \right)^2
  =
  C  \zeta^{2 + 2H}.
\]

It follows immediately from   \cref{lem:negative moments} and   \cref{negative moments-1} that

\begin{equation*}
   \mathrm{E} \left( ||D^B  \sigma ^2 _{_Y}|| _{H} \right) ^{-8}\leq C \mathrm{E} \zeta^{-8 - 8H}
 \leq C \mathrm{E} \left( \tau _{\varepsilon} \wedge \varepsilon\wedge\frac{T}{2}
 \right) ^{-8 - 8H}<\infty.
\end{equation*}

 \qed

\textbf{Proof of   \cref{thm:density fBm}}.

 From   \cref{lem:lurid} and Proposition 2.1.1, \cite{Nua06} we get the first part of equality \cref{eq:loin-1}:
 \begin{equation*}
  p_{\sigma ^2 _{_Y}} (u) = \mathrm{E} \left[ \I_{\sigma ^2 _{_Y} >u  }
  \delta \left(
  \frac{D^B \sigma  ^2 _{_Y}}{||D  \sigma ^2 _{_Y}||^2 _{H}}\right)
  \right] .
 \end{equation*}

 To get the second part, note that $\eta:=\left(||D  \sigma ^2 _{_Y}|| _{H}\right)^{-2}$ admits stochastic derivative and, according to Proposition 1.3.3 from \cite{Nua06},
 the following holds
  \begin{equation*}\begin{gathered}
  \delta \left(
  \frac{D^B \sigma  ^2 _{_Y}}{||D^B  \sigma ^2 _{_Y}||^2 _{H}}\right)=\int_0^T\eta D^B_u (\sigma  ^2 _{_Y})dB_u=\eta \int_0^T  D^B_u (\sigma  ^2 _{_Y})dB_u\\-
  \int_0^TD^B_u \eta  D^B_u (\sigma  ^2 _{_Y})du=2 \eta\int_0^T u ^{1/2 - H}\int\limits _u ^T \kappa(Y _s)
 l(u,s)dsdB_u-
  \int_0^TD^B_u \eta  D^B_u (\sigma  ^2 _{_Y})du.\end{gathered}\end{equation*}
 According to Lemma 2.10 from \cite{leon}, we can apply the Fubini theorem for the Skorokhod integral.  Then  $$\int_0^T u ^{1/2 - H}\int\limits _u ^T \kappa(Y _s)
 l(u,s)dsdB_u=\int_0^T\kappa(Y _s)(\int_0^su ^{1/2 - H}l(u,s)dB_u)ds,$$
 where the interior integral is a Wiener one.

 Finally, taking into account that  $m _{_Y} =  X_0 + bT- \frac 12 \sigma ^2 _{_Y}$, we get

 \[
  \mathrm{E}
 \frac { (x - m _{_Y}) }{\sigma  _{_Y} ^3 \sqrt{2\pi}}
        \exp  \Bigg\{
    -\frac{ (x - m _{_Y}) ^2}{\sigma _{_Y} ^2}
       \Bigg\}
       \]
       \[=
       \int_ {\R}
       \frac { (x + u/2 - X_0 - bT) }{u ^3 \sqrt{2\pi}}
        \exp  \Bigg\{
    -\frac{ (x + u/2 - X_0 - bT) ^2}{u ^2}
       \Bigg\} p_{\sigma ^2 _{_Y}} (u) du.
 \]

 Combining this with \cref{eq:assuage}, we get the proof.

\section{Simulations}\label{sims}

 In this section we use the discretization schemes
 proposed in \cref{1st level} and \cref{2nd level} to simulate
the option price. We treat  double and single discretization, respectively.

The  values of $b, \alpha$ and $T$
are the same in all simulations, and equal   $b = 0.2, \alpha = 0.6, T = 1$.
In \cref{tab:11} and \cref{tab:22} we give the results of simulations
based on \cref{eq:for f} (double discretization)
for different $n$ and $\sigma$. The functions $f$
and the values of $H$ are given in the table headers.

\begin{table}[h]
\caption{double discretization,  $f(s) = (s-1)_+ +  \I_{s > 1 }$, $H =0.6$ }
\label{tab:11}
\centering
\begin{tabular}{|c|c|c|c|c|c|c|c|} \hline
$n$                              & 125      & 250     & 500      & 1000      & 2000      & 4000    & 8000 \\ \hline

$\sigma(y) = \sqrt{|y| + 0.1}$   & 0.95124  & 0.92149 & 0.93664  & 0.89628   & 0.88124   & 0.89717 & 0.92390    \\

$\sigma(y) = |y|+0.1$            & 0.92121  & 0.95820 & 0.94733  & 1.02572   & 0.90530   & 0.92062 & 0.97430    \\
$\sigma(y) = \sqrt{y^2 + 1}$     & 0.93357  & 1.01340 & 0.99205  & 0.95801   & 0.97705   & 0.96882 & 0.92312   \\
$\sigma(y) = \sin ^2 (y) $       & 0.87957  & 0.87842 & 0.94525  & 0.93053   & 0.91097   & 0.89368 & 1.00256  \\
$+ 0.05$                         &          &         &           &          &          &         &          \\
\hline
\end{tabular}
\end{table}

\begin{table}[h]
\caption{double discretization,  $f(s) = (s-1.5)_+ +  \I_{s > 2 }$, $H =0.8$ }
\label{tab:22}
\centering
\begin{tabular}{|c|c|c|c|c|c|c|c|} \hline
$n$                              & 125      & 250     & 500      & 1000      & 2000      & 4000    & 8000 \\ \hline
$\sigma(y) = \sqrt{|y| + 0.1}$   & 0.48073  & 0.46643 & 0.53185  & 0.53124   & 0.53128   & 0.50020 & 0.57804    \\
$\sigma(y) = |y|+0.1$            & 0.51931  & 0.48957 & 0.50875  & 0.48999   & 0.46509   & 0.48525 & 0.55298    \\
$\sigma(y) = \sqrt{y^2 + 1}$     & 0.66845  & 0.81584 & 0.66681  & 0.64368   & 0.65982   & 0.76746 & 0.74945   \\
$\sigma(y) = \sin ^2 (y) $       & 0.31951  & 0.30628 & 0.28770  & 0.29487   & 0.30395   & 0.29509 & 0.30490  \\
$+ 0.05$                         &          &         &           &          &          &         &          \\
\hline
\end{tabular}
\end{table}

 \cref{tab:33} and \cref{tab:44} present the results of simulations
 for the same parameters as in \cref{tab:11} and \cref{tab:22}, respectively,
 but the price is computed using \eqref{eq:assuage},
 which corresponds to single discretization.

\begin{table}[h]
\caption{single discretization,  $f(s) = (s-1)_+ +  \I_{s > 1 }$, $H =0.6$ }
\label{tab:33}
\centering
\begin{tabular}{|c|c|c|c|c|c|c|c|} \hline
$n$                              & 125      & 250     & 500      & 1000      & 2000      & 4000    & 8000 \\ \hline

$\sigma(y) = \sqrt{|y| + 0.1}$   & 0.92266  & 0.92257 & 0.92249  & 0.92232   & 0.92219   & 0.92263 & 0.92262    \\

$\sigma(y) = |y|+0.1$            & 0.92985  & 0.92910 & 0.92980  & 0.92937   & 0.92955   & 0.92958 & 0.92916    \\
$\sigma(y) = \sqrt{y^2 + 1}$     & 0.96127  & 0.96115 & 0.96124  & 0.96114   & 0.96155   & 0.96107 & 0.96093  \\
$\sigma(y) = \sin ^2 (y) $       & 0.92171  & 0.92143 & 0.92148  & 0.92136   & 0.92086   & 0.92092 & 0.92131  \\
$+ 0.05$                         &          &         &           &          &          &         &          \\
\hline
\end{tabular}
\end{table}

\begin{table}[h]
\caption{single discretization,  $f(s) = (s-1.5)_+ +  \I_{s > 2 }$, $H =0.8$ }
\label{tab:44}
\centering
\begin{tabular}{|c|c|c|c|c|c|c|c|} \hline
$n$                              & 125      & 250     & 500      & 1000      & 2000      & 4000    & 8000 \\ \hline

$\sigma(y) = \sqrt{|y| + 0.1}$   & 0.51543  & 0.51558 & 0.51415  & 0.51461   & 0.51592   & 0.51498 & 0.51525   \\

$\sigma(y) = |y|+0.1$            & 0.50348  & 0.50743 & 0.50659  & 0.50460   & 0.50591   & 0.50450 & 0.50641    \\
$\sigma(y) = \sqrt{y^2 + 1}$     & 0.67837  & 0.67881 & 0.67782  & 0.67843   & 0.67846   & 0.67848 & 0.67825   \\
$\sigma(y) = \sin ^2 (y) $       & 0.31314  & 0.31239 & 0.30948  & 0.31214   & 0.30948   & 0.30658 & 0.30887  \\
$+ 0.05$                         &          &         &           &          &          &         &          \\
\hline
\end{tabular}
\end{table}

Note that the results of simulations related to a single discretization  appear to be
much more consistent in $n$. This is maybe
 due to the fact that, unlike in \cref{eq:for f}, the value
under the sign of expectation in \cref{eq:assuage} is bounded,
and thus the average over $10^4$ trials
gives a good approximation.
 The "geometric" nature of  $S$ is probably also a factor. Indeed, let us
have a look at the $10$ largest values of $S_T$ in $10^4$ trials
 with the $\sigma(y) = |y|+0.1$, $H =0.8$, and $n = 400$:
\[
89.4301 \ \  36.2412   \ \ 34.7761 \ \  34.7639 \ \  34.7330 \ \
30.7971 \ \  27.0400  \ \  25.7752 \ \  24.2836 \ \  23.0231.
\]
Taking into account that $F$ exhibits quadratic growth,
we see that the result of a single trial can influence
the average over all $10^4$ trials.
The consistency of single discretization  is also confirmed by the results
in the \cref{tab:55}, where we also represent the results of simulation of  the direct
average   $\mathrm{E}f(S_T)$, without reducing to continuous functions,  over the same number of trials.
Note that we take the same realizations for double discretization,
single discretization and the direct average, observing
the evident  correlation between the results for double discretization  and for the direct average.

\begin{table}[h]
\caption{$f(s) = (s-1)_+ +  \I_{s > 1 }$, $\sigma(y) = |y|+0.2$, $H =0.75$ }
\label{tab:55}
\centering
\begin{tabular}{|c|c|c|c|c|c|c|c|} \hline
$n$             & 125      & 250     & 500      & 1000      & 2000      & 4000    & 8000 \\ \hline
level one       & 0.96675   & 0.90599  & 0.90634  & 0.9483  & 1.0206  & 0.91574  & 0.92130   \\
level two     & 0.93501  & 0.93429  & 0.93415  & 0.93449 &  0.93471  & 0.93402 &  0.93415   \\
direct average & 0.9602  & 0.9047  & 0.91156  & 0.94228 &  1.0023 &  0.92594  & 0.91595    \\
\hline
\end{tabular}
\end{table}

We thus conjecture that in practice single discretization  gives a good approximation
for the expectation $\mathrm{E}f(S_T)$,
better than double discretization  or the direct averaging.
The simulations are provided  in MATLAB.
To simulate the fBm, we use the algorithm
"Fractional Brownian motion generator"
by Zdravko Botev
which can be found at \cite{fbm1d}.
For simulations related to double discretization  we take the average
of the value  under the expectation in the right hand side
of \cref{eq:for f}
over $10^4$ trials. For simulations related to single discretization,
we replace infinite interval of integration in the right hand side
of \eqref{eq:assuage} with a finite one, making sure
that the integral over the complement is small. Then
we discretize the finite interval; the partition size is $2500$.
For each $x$ from the partition, we take
the average over $10^4$ trials of the value
under the expectation in \eqref{eq:assuage}.
The trials are common for all $x$, i.e.
it is not necessary to generate $10^4$
trials for every $x$ from the partition.

\bibliographystyle{abbrv}
\bibliography{Citations}

\end{document}

%% file: title_page.tex
\usepackage{lipsum}
\usepackage{amsfonts}
\usepackage{graphicx}
\usepackage{epstopdf}
\usepackage{algorithmic}
\ifpdf
  \DeclareGraphicsExtensions{.eps,.pdf,.png,.jpg}
\else
  \DeclareGraphicsExtensions{.eps}
\fi

\newcommand{\TheTitle}{Option pricing with fractional stochastic volatility and discontinuous payoff function of polynomial growth}
\newcommand{\TheAuthors}{V. Bezborodov, L. Di Persio, and Y. Mishura}

 \headers{Option pricing with fractional stochastic volatility}{\TheAuthors}

\title{{\TheTitle}\thanks{21 July 2016
\funding{The research of Viktor Bezborodov was  supported by
the Department of Computer Science of the University of Verona,
  Luca Di Persio was supported by the Gruppo Nazionale per
l'Analisi Matematica, la Probabilit\`{a} e le loro Applicazioni (GNAMPA)   within the project
"Stochastic Partial Differential Equations and Stochastic Optimal Transport
with Applications to Mathematical Finance", and
Yuliya Mishura is thankful to the University of Verona
for the hospitality during her visit.}}}

\author{
  Viktor Bezborodov\thanks{The University of Verona,
  The Department of Computer Science
    (\email{viktor.bezborodov@univr.it})}
      \and
  Luca Di Persio \thanks{The University of Verona,
  The Department of Computer Science
    (\email{luca.dipersio@univr.it})}
  \and
  Yuliya Mishura\thanks{Taras Shevchenko National University of Kyiv,
  The Faculty of Mechanics and Mathematics
  (\email{myus@univ.kiev.ua})}
}

\usepackage{amsopn}
